\documentclass[12pt,reqno]{amsart}
\usepackage{amscd}
\usepackage{amssymb}
\usepackage[mathscr]{eucal}
\usepackage{enumerate}

\makeatletter

\def\subsection{\@startsection{subsection}{2}%
  \z@{.5\linespacing\@plus.7\linespacing}{1pt}%
      {\normalfont\bfseries}}

\def\l@section{\@tocline{1}{0pt}{1pc}{}{\bfseries}}
\def\l@subsection{\@tocline{2}{0pt}{3.1em}{5pc}{}}

\makeatother

\allowdisplaybreaks[1]

\newtheorem{thm}{Theorem}[section]
\newtheorem{lem}[thm]{Lemma}
\newtheorem{cor}[thm]{Corollary}

\theoremstyle{definition}

\newtheorem{rem}[thm]{Remark}
\newtheorem{defn}[thm]{Definition}

\theoremstyle{remark}

\numberwithin{equation}{section}

\newcommand{\C}{{\mathbb{C}}}
\newcommand{\R}{{\mathbb{R}}}

\newcommand{\N}{{\mathbb{N}}}
\newcommand{\T}{{\mathbb{T}}}

\DeclareMathOperator{\Ad}{Ad}
\DeclareMathOperator{\Tr}{Tr}

\DeclareMathOperator{\id}{id}

\DeclareMathOperator{\spa}{span}

\DeclareMathOperator{\Irr}{Irr}

\DeclareMathOperator{\Ima}{Im}

\def\mD{\mathcal{D}}

\def\meQ{\mathscr{Q}}

\def\meV{\mathscr{V}}

\def\mfg{\mathfrak{g}}

\def\al{\alpha}
\def\be{\beta}
\def\ga{\gamma}
\def\de{\delta}
\def\ka{\kappa}
\def\la{\lambda}

\def\vep{\varepsilon}
\def\ph{{\phi}}
\def\ps{{\psi}}
\def\vph{\varphi}
\def\vrho{\varrho}
\def\om{\omega}
\def\si{\sigma}
\def\ta{\tau}
\def\th{\theta}

\def\Th{\Theta}

\def\De{\Delta}
\def\vDe{\varDelta}
\def\La{\Lambda}

\def\el{\ell}

\def\ovl{\overline}                 
\def\wdh{\widehat}

\def\op{{\mathrm{op}}}

\def\hB{{\widehat{B}}}
\def\hC{{\widehat{C}}}

\def\hJ{\hat{J}}

\def\hS{\hat{S}}
\def\hR{\hat{R}}

\def\hta{{\hat{\ta}}}

\def\hvep{{\hat{\vep}}}

\def\opi{{\ovl{\pi}}}                

\def\subs{\subset}                   

\def\setm{\setminus}

\def\oti{\otimes}                    
\def\rti{\rtimes}                    

\def\col{\colon}
\def\ra{\hspace{-0.5mm}\rightarrow\!}

\def\btr{{\boldsymbol{1}}}                   
\def\bG{\mathbb{G}}
\def\bH{\mathbb{H}}
\def\bK{\mathbb{K}}

\def\bhG{{\wdh{\mathbb{G}}}}

\def\ltG{L^2(\bG)}
\def\AG{A(\bG)}
\def\CG{C(\bG)}
\def\CGq{C(\bG_q)}
\def\lG{L^\infty(\bG)}                      
\def\LG{L(\bG)}                             
\def\RG{R(\bG)}                             

\def\lH{L^\infty(\bH)}
\def\LH{L(\bH)}
\def\LK{L(\bK)}
\def\RH{R(\bH)}
\def\AH{A(\bH)}
\def\CH{C(\bH)}
\def\CK{C(\bK)}

\def\IG{\Irr(\bG)}
\def\IH{\Irr(\bH)}

\begin{document}

\title{A Characterization of right coideals of quotient type 
and its application to classification of Poisson boundaries}

\author{Reiji Tomatsu}
\address{Reiji Tomatsu, Department of Mathematical Sciences,
University of Tokyo, 3-8-1 Komaba, Meguro, Tokyo, 153-8914, JAPAN}
\email{tomatsu@ms.u-tokyo.ac.jp}


\maketitle

\begin{abstract}
Let $\bG$ be a co-amenable compact quantum group. 
We show that 
a right coideal of $\bG$ is of quotient type 
if and only if 
it is the range of a conditional expectation preserving the Haar state 
and is globally invariant under 
the left action of the dual discrete quantum group. 
We apply this result to theory of Poisson boundaries 
introduced by Izumi for discrete quantum groups 
and generalize a work of Izumi-Neshveyev-Tuset on $SU_q(N)$ 
for co-amenable compact quantum groups 
with the commutative fusion rules. 
More precisely, 
we prove that the Poisson integral is an isomorphism 
between the Poisson boundary 
and the right coideal of quotient type 
by maximal quantum subgroup of Kac type. 
In particular, 
the Poisson boundary and the quantum flag manifold 
are isomorphic for any $q$-deformed classical compact Lie group. 
\end{abstract}

\section{Introduction}

Since Woronowicz introduced the axiomatic 
compact quantum groups \cite{Wor3}, 
they have attracted a growing interest of many researchers 
as a framework to describe new types of symmetries. 
In this paper, we study two typical examples of ergodic actions of 
compact quantum groups, namely right coideals and Poisson boundaries. 

Let $\bG$ be a compact quantum group. 
A right coideal is a von Neumann subalgebra of the function algebra 
on $\bG$ which is globally invariant under 
the right translation action of $\bG$. 
Taking the fixed point algebra 
of a left action of a quantum subgroup gives 
an example of a right coideal. 
We say that such a right coideal is of quotient type. 
When $\bG$ is an ordinary group, 
it is well-known that all the right coideals are of quotient type. 
However, when $\bG$ is a quantum group, 
not all the right coideals 
are realized as quotients \cite{Po1}, \cite{Po2}, \cite{To2}. 
This fact presents a contrast between the quantum groups
and the ordinary ones. 

Another difference 
occurs in behaviors of infinite tensor product actions. 
In the ordinary case, such an action is minimal 
and, in particular, 
the relative commutant of the fixed point algebra is trivial. 
However, this is not the case for quantum groups. 
In \cite{I1}, Izumi has described this contrast 
by introducing the notion of a Poisson boundary 
of a dual discrete quantum group. 
More precisely, 
he has shown that the Poisson boundary is isomorphic 
to the relative commutant of the fixed point algebra. 
Moreover he has also studied the Poisson boundary of 
the dual of $SU_q(2)$ 
and shown a striking result 
that the Poisson boundary is isomorphic 
to the standard Podle\'{s} sphere $L^\infty(\T\setm SU_q(2))$ 
introduced in \cite{Po1}, \cite{Po2}. 
On one hand, this result has led to a conjecture
that the Poisson boundary 
is isomorphic to the quantum flag manifold 
for any $q$-deformed classical compact Lie group. 
For $SU_q(N)$, 
the conjecture was confirmed affirmatively 
in \cite{INT}. 
Note that a $q$-deformed classical compact Lie group is 
co-amenable, that is, the dual discrete quantum group 
is amenable in the sense of \cite{BeMuTu1}. 
On the other hand for non-amenable cases, 
the Poisson and Martin boundaries \cite{NT} of 
universal orthogonal discrete quantum groups are studied 
in \cite{VV} and \cite{VVe}. 

In this paper, 
we first characterize 
when a right coideal is of quotient type 
and 
second apply the characterization 
to determine Poisson boundaries of 
amenable discrete quantum groups with the commutative fusion rules. 
A right coideal of quotient type has the following 
two properties. 
The first one is the \textit{expectation property}, 
namely, existence of a normal conditional expectation 
preserving the Haar state from the function algebra. 
The second one is the \textit{coaction symmetry} 
which means that the left action of the dual preserves 
the right coideal. 
Assuming amenability of the dual, 
we can prove the following theorem (Theorem \ref{thm: exco}). 

\vspace{10pt}
\hspace{-12.2pt}\textbf{Theorem 1.} 
Let $\bG$ be a co-amenable compact quantum group 
and $B\subs \lG$ a right coideal. 
Then 
$B$ is of quotient type if and only if 
$B$ has the expectation property and the coaction symmetry. 

\vspace{10pt}
Next we study a Poisson boundary of an amenable discrete quantum group. 
In order to compute a Poisson boundary, 
we present an approach which differs from that of \cite{INT}. 
The key point of our proof is 
to construct an ``inverse''of the Poisson integral. 
Although this strategy is the same as the one taken in \cite{INT}, 
we do it by utilizing not the Berezin transforms 
but an invariant mean of a dual discrete quantum group. 
Then using Theorem 1, 
we show that the Poisson boundary is isomorphic 
to a right coideal of quotient type by a quantum subgroup. 
Moreover, we can specify the quantum subgroup which is the 
\textit{maximal} quantum subgroup of Kac type with respect
to inclusions. 
After this work was done, we learned from S. Vaes that 
the notion 
has been already introduced as the canonical Kac quotient in \cite{So}. 
Our main result is the following theorem (Theorem~\ref{thm: poisson}).  

\vspace{10pt}
\hspace{-12.2pt}\textbf{Theorem 2.} 
Let $\bG$ be a co-amenable compact quantum group. 
Assume that its fusion algebra is commutative. 
Then the following statements hold. 
\vspace{5pt}

(1) There exists a unique maximal quantum subgroup of Kac type $\bH$. 

(2) The Poisson integral 
$\Th\col L^\infty(\bH\setm\bG)\ra H^\infty(\bhG)$ 
is an isomorphism. 

\vspace{10pt}
In particular, $\bG$ is of Kac type if and only if $H^\infty(\bhG)=\C$. 
This yields the minimality of an infinite tensor product actions of $\bG$. 
For a $q$-deformed classical compact Lie group $\bG_q$, 
the maximal quantum subgroup of Kac type is the maximal torus $\T$. 
Therefore, we obtain the following result (Corollary \ref{cor: flag}). 

\vspace{10pt}
\hspace{-12.2pt}\textbf{Theorem 3.} 
Let $\bG_q$ be the $q$-deformation of a classical compact Lie group $\bG$. 
Then the Poisson integral 
$\Th\col L^\infty(\T\setm\bG_q)\ra H^\infty(\widehat{\bG_q})$ 
is an isomorphism. 

\vspace{10pt}
\hspace{-12.15pt}\textbf{Notations.} 
Let $M$ be a von Neumann algebra with predual $M_*$. 
We denote by $M_*^+$ the set of positive functionals in $M_*$. 
For a linear functional $\th$ on $M$, we define a linear functional 
$\ovl{\th}$ by $\ovl{\th}(x)=\ovl{\th(x^*)}$ 
for $x\in M$. 
For a weight $\ph$ on $M$, 
we set $n_\ph=\{x\in M\mid \ph(x^* x)<\infty\}$, 
$m_\ph=n_\ph^* n_\ph$. 
We denote by $m_\ph^+$ the set of positive elements in $m_\ph$. 
For a linear subspace $X\subs M$, we denote by $\ovl{X}^{\,\rm{w}}$ 
the weak closure of $X$. 
We denote by $\oti$ the minimal tensor product for $C^*$-algebras 
and the spatial tensor product for von Neumann algebras. 

\section{Preliminaries}

We collect necessary results on compact (discrete) quantum groups. 

\subsection{Compact quantum groups.}\label{subsect: cqg}
Our standard references are 
\cite{BS}, \cite{I1}, \cite{Wor3}. 
For the notion of a compact quantum group, 
we adopt the definition in \cite[Definition 2.1]{Wor3} as follows: 

\begin{defn}
A \textit{compact quantum group} $\bG$ is a pair 
$(C(\bG),\delta_\bG)$ which satisfies the following conditions: 

\begin{enumerate}

\item $C(\bG)$ is a separable unital $C^*$-algebra. 

\item 
The map $\delta_\bG\col C(\bG)\ra C(\bG)\oti C(\bG)$ 
is a coproduct, i.e. it is a faithful unital $*$-homomorphism satisfying 
the coassociativity condition, 
\[(\delta_\bG\otimes\id)\circ \delta_\bG
=(\id\otimes\delta_\bG)\circ\delta_\bG.
\]

\item 
The vector spaces $\delta_\bG(\CG)(\C\otimes \CG)$ 
and $\delta_\bG(\CG)(\CG\otimes\C)$ 
are dense in $\CG\otimes \CG$.
\end{enumerate}
\end{defn}

Let $h_\bG$ be
the \textit{Haar state} on $\CG$ which satisfies the invariance 
condition, 
\[
(\id\oti h_\bG)(\de_\bG(a))=h_\bG(a)1=(h_\bG\oti\id)(\de_\bG(a))
\quad
\mbox{for all}\ a\in \CG. 
\]
In this paper, we always assume that the Haar states are faithful. 
If the Haar state is tracial, we say that the compact quantum group 
is of \textit{Kac type} \cite{EnSch}. 
Let $(\pi_h, \ltG,\hat{1}_h)$ be the GNS triple of $h_\bG$, 
which consists of the representation, the Hilbert space 
and the GNS cyclic vector, respectively. 
We always omit $\pi_h$ and regard $\CG$ as a $C^*$-subalgebra of $B(\ltG)$. 
We set a von Neumann algebra $\lG=\ovl{\CG}^{\,\rm{w}}$. 

The \emph{multiplicative unitaries} $V_\bG$ and $W_\bG$ are defined by 
\begin{equation}\label{eq: W}
W_\bG^*(x\hat{1}_h\oti y\hat{1}_h )=\de_\bG(y)(x \hat{1}_h \oti \hat{1}_h)
\quad \mbox{for}\ x,y\in \CG, 
\end{equation}
\begin{equation}\label{eq: V}
V_\bG(x\hat{1}_h\oti y\hat{1}_h )=\de_\bG(x)(\hat{1}_h \oti y\hat{1}_h)
\quad \mbox{for}\ x,y\in \CG. 
\end{equation}
Then we have the pentagon equalities, 
\begin{equation}\label{eq: Wpenta}
(W_\bG)_{12}(W_\bG)_{13}(W_\bG)_{23}=(W_\bG)_{23}(W_\bG)_{12},
\end{equation}
\begin{equation}\label{eq: Vpenta}
(V_\bG)_{12}(V_\bG)_{13}(V_\bG)_{23}=(V_\bG)_{23}(V_\bG)_{12}. 
\end{equation}
Using them, 
we can extend the coproduct $\de_\bG$ to $\lG$ by 
\[
\de_\bG(x)=V_\bG(x\oti1)V_\bG^*=W_\bG^*(1\oti x)W_\bG
\quad\mbox{for}\  x\in \lG. 
\] 
The normal state $h_\bG(x)=(x\hat{1}_h,\hat{1}_h)$, $x\in\lG$ 
is invariant under $\de_\bG$. 
Since $h_\bG$ on $\CG$ has the modular automorphism group 
\cite[Theorem 2.4]{Wor3}, 
$h_\bG$ is also faithful on $\lG$. 
Then the pair $(\lG,\de_\bG)$ is a 
\emph{von Neumann algebraic compact quantum group} 
in the sense of \cite{KV}. 

Let $H$ be a Hilbert space and $v\in B(H)\oti \lG$ a unitary. 
If $(\id\oti\de_\bG)(v)=v_{12}v_{13}$, we say that $v$ is a 
right \emph{unitary representation} of $\bG$. 
A left unitary representation is similarly defined. 
For example, 
the unitaries $W_\bG$ and $V_\bG$ are left and right unitary 
representations, respectively. 
Let $v\in B(H)\oti \lG$ be a unitary representation. 
Suppose that $T\in B(H)$ satisfies $(T\oti1)v=v(T\oti1)$. 
If any such element $T$ must be a scalar, 
$v$ is said to be \emph{irreducible}. 
Any unitary representation is completely decomposable, 
that is, it is a direct sum of irreducible ones. 
The set of the equivalence classes of all the irreducible representations 
is denoted by $\IG$. 
For $\pi\in \IG$, we choose 
a representation Hilbert space $H_\pi$ 
and 
an irreducible representation 
$v_\pi=(v_{\pi_{i,j}})_{i,j\in I_\pi}\in B(H_\pi)\oti \CG$. 
We call an irreducible representation $1\in \C\oti \lG$ 
the \emph{trivial representation} 
and denote by $\btr$ 
the equivalence class. 

We define a dense unital $*$-subalgebra $\AG\subs\CG$ 
by 
\[
\AG=
\spa\{v_{\pi_{i,j}}\mid i,j\in I_\pi,\ \pi\in\IG\}. 
\] 
It is clear that $\de_\bG(\AG)\subs \AG\oti \AG$. 
We define the Hopf algebra structure, namely, 
the \emph{antipode} $\ka_\bG$ and the \emph{counit} $\vep_\bG$ 
on $\AG$ as follows. 
The invertible antimultiplicative map $\ka_\bG\col \AG\ra\AG$ 
is defined by 
\begin{equation}\label{eq: antipode}
\ka_\bG(v_{\pi_{i,j}})=v_{\pi_{j,i}}^*
\quad\mbox{for}\ 
i,j\in I_\pi,\ \pi\in\IG. 
\end{equation}
The unital $*$-homomorphism $\vep_\bG\col \AG\ra \C$ 
is defined by 
\[
\vep_\bG(v_{\pi_{i,j}})=\de_{i,j}
\quad\mbox{for}\ 
i,j\in I_\pi,\ \pi\in\IG. 
\]
In fact $\AG$ is a Hopf $*$-algebra, that is, 
$\ka_\bG(\ka_\bG(x)^*)^*=x$ holds for any $x\in\AG$. 
For any finite dimensional unitary representation $v\in B(H)\oti \AG$, 
we have 
$(\id\oti\ka_\bG)(v)=v^*$ and $(\id\oti\vep_\bG)(v)=1$, 
which follow from the complete decomposability of $v$. 

We introduce the \emph{Woronowicz characters} 
$\{f_z^\bG\}_{z\in\C}$ on $\AG$ \cite[Theorem2.4]{Wor3}. 
The multiplicative functional $f_z^\bG\col \AG\ra\C$ 
satisfies the following properties: 

\begin{enumerate}[(i)]

\item $f_0^\bG=\vep_\bG$. 

\item For any $a\in \AG$, 
the function $\C\ni z \mapsto f_z^\bG(a)\in\C$ 
is entirely holomorphic.

\item 
$(f_{z_{1}}^\bG\oti f_{z_{2}}^\bG)\circ\de_\bG=f_{z_{1}+z_{2}}^\bG$ 
for any $z_{1}, z_{2}\in\C$. 

\item For any $z\in \C$ and $a\in \AG$, 
$f_z^\bG(\ka(a))=f_{-z}^\bG(a)$, 
$f_{z}^\bG(a^*)=\ovl{f_{-\bar{z}}^\bG(a)}$.

\item For any $a\in \AG$, 
$\ka_\bG^2(a)=(f_1^\bG\oti\id\oti f_{-1}^\bG)
\big{(}(\de_\bG\oti\id)(\de_\bG(a))\big{)}$.

\item 
For any $a,b\in\AG$, 
$h_\bG(ab)=h_\bG\big{(}b\,(f_1^\bG\oti\id\oti f_{1}^\bG)
\big{(}(\de_\bG\oti\id)(\de_\bG(a))\big{)}\big{)}$. 
\end{enumerate}

The modular automorphism group $\{\si_t^{h_\bG}\}_{t\in\R}$ is given by 
\[
\si_t^{h_\bG}(x)
=(f_{it}^\bG\oti\id\oti f_{it}^\bG)\big{(}(\de_\bG\oti\id)(\de_\bG(x))\big{)}
\quad\mbox{for all}\ 
t\in \R,\ x\in \AG. 
\]
We define the following map $\ta_t^\bG\col \AG\ra\AG$ by
\[
\ta_t^\bG(x)
=(f_{it}^\bG\oti\id\oti f_{-it}^\bG)\big{(}
(\de_\bG\oti\id)(\de_\bG(x))\big{)}
\quad\mbox{for all}\ 
t\in \R,\ x\in \AG. 
\]
Then $\{\ta_t^\bG\}_{t\in\R}$ is a one-parameter automorphism group on $\AG$ 
and it is called the \emph{scaling automorphism group}. 
Note that any element of $\AG$ is analytic for $\{\si_t^{h_\bG}\}_{t\in\R}$ 
and $\{\ta_t^\bG\}_{t\in\R}$. 
Define a map $R_\bG\col \AG\ra\AG$ by 
\[
R_\bG(x)=\ka_\bG(\ta_{i/2}^\bG(x))
\quad\mbox{for}\ 
x\in \AG. 
\]
Then $R_\bG$ is a $*$-antiautomorphism on $\AG$ with $R_\bG^2=\id$. 
We call $R_\bG$ the \emph{unitary antipode}. 
Actually $R_\bG$ commutes with $\ta_t^\bG$ for all $t\in\R$, 
and it satisfies
\[
\ka_\bG=R_\bG\circ \ta_{-i/2}^\bG=\ta_{-i/2}^\bG\circ R_\bG. 
\]
Since the Haar state $h_\bG$ is invariant 
under 
the $*$-preserving maps $\si_t^{h_\bG}$, $\ta_t^\bG$ and $R_\bG$, 
we can extend them to the maps on $\CG$, moreover on $\lG$. 
Then the following relations among $\de_\bG$ and them hold. 
\begin{equation}\label{eq: de1}
h_\bG\circ\ta_t^\bG=h_\bG=h_\bG\circ R_\bG, \quad 
\de_\bG\circ\si_t^{h_\bG}=(\si_t^{h_\bG}\oti\ta_{-t}^\bG)\circ\de_\bG, 
\end{equation}
\begin{equation}\label{eq: de2}
\de_\bG\circ \ta_t^\bG=(\ta_t^\bG\oti\ta_t^\bG)\circ\de_\bG, 
\quad 
\de_\bG\circ R_\bG=(R_\bG\oti R_\bG)\circ \de_\bG^\op, 
\end{equation}
where $\de_\bG^\op(x)=\de_\bG(x)_{21}$. 

Let $v\in B(H)\oti \AG$ be a finite dimensional unitary representation. 
Set $F_v^\bG=(\id\oti f_1^\bG)(v)$. 
Then $F_v^\bG$ satisfies the following properties: 
\begin{enumerate}[(i)]

\item $F_v^\bG$ is a non-singular positive operator on $H$. 

\item For any $z\in\C$, $(\id\oti f_z^\bG)(v)=(F_v^\bG)^z$. 

\item For any $t\in\R$, 
$(\id\oti\ta_t^\bG)(v)=((F_v^\bG)^{it}\oti1)v((F_v^\bG)^{-it}\oti1)$. 
\item 
For any $t\in\R$, 
$(\id\oti\si_t^{h_\bG})(v)=((F_v^\bG)^{it}\oti1)v((F_v^\bG)^{it}\oti1)$. 

\item 
Let $w\in B(H_w)\oti\AG$ be a finite dimensional 
unitary representation of $\bG$. 
If a linear map $T\col H_v\ra H_w$ satisfies $(T\oti1)v=w(T\oti1)$, 
then $T F_v=F_w T$. 

\end{enumerate}

For $\pi\in\IG$, we write $F_\pi^\bG$ (or simply $F_\pi$) 
for $F_{v_\pi}^\bG$. 
Let $\Tr_\pi$ be the non-normalized trace on $B(H_\pi)$. 
Set $D_\pi=\Tr_\pi(F_\pi)$. 
Then we have the orthogonal relations, 
\[
h_\bG(v_{\pi_{i,j}} v_{\rho_{k,\el}}^*)
=D_{\pi}^{-1}F_{\pi_{j,j}}\delta_{\pi,\rho}\delta_{i,k}\de_{j,\el},
\quad
h_\bG(v_{\pi_{i,j}}^{*}v_{\rho_{k,\el}})
=D_{\pi}^{-1}F_{\pi_{i,i}}^{-1}\delta_{\pi,\rho}\de_{i,k}\delta_{j,\el}. 
\]

We decompose $W_\bG$ and $V_\bG$ 
into irreducible representations.  
For $\pi\in\IG$, we define two systems of matrix units 
$\{e_{\pi_{i,j}}\}_{i,j\in I_\pi}$ and $\{f_{\pi_{i,j}}\}_{i,j\in I_\pi}$ 
in $B(\ltG)$ by 
\begin{equation}\label{eq: epi}
e_{\pi_{i,j}} (v_{\rho_{k,\el}}\hat{1}_h)
=\de_{\pi,\rho}\de_{j,\el}\, 
v_{\rho_{k,i}}\hat{1}_h
\quad\mbox{for}\ k,\el\in I_\rho,\ \rho\in\IG,
\end{equation}
\begin{equation}\label{eq: fpi}
f_{\pi_{i,j}} (v_{\rho_{k,\el}}^*\hat{1}_h)
=\de_{\pi,\rho}\de_{j,k}\, 
v_{\rho_{i,\el}}^*\hat{1}_h 
\quad\mbox{for}\ k,\el\in I_\rho,\ \rho\in\IG. 
\end{equation}
Then we have 
\begin{equation}\label{eq: decomp}
W_\bG=\sum_{\pi\in\IG}\sum_{i,j\in I_\pi} 
v_{\pi_{i,j}}\oti f_{\pi_{i,j}},
\quad 
V_\bG=\sum_{\pi\in\IG}\sum_{i,j\in I_\pi}
e_{\pi_{i,j}}\oti v_{\pi_{i,j}}. 
\end{equation}
Setting $\pi=\btr$ at (\ref{eq: epi}), we have 
\[
e_\btr (x \hat{1}_h)= h_\bG(x)\hat{1}_h 
\quad\mbox{for}\ x\in \lG. 
\]
The projection $e_\btr=f_\btr$ is minimal in $B(\ltG)$. 

We will need some relations among $W_\bG$, $V_\bG$ and the modular 
objects of $\lG$. 
Let $\vDe_{h_\bG}$ and $J_\bG$ be the modular operator 
and the modular conjugation of $h_\bG$. 
We define a conjugate unitary $\hJ_\bG$ by 
$\hJ_\bG(x\hat{1}_h)=R_\bG(x^*)\hat{1}_h$ for $x\in \CG$. 
We set a unitary $U_\bG=J_\bG\hJ_\bG=\hJ_\bG J_\bG$. 
The implementing unitary $P_\bG^{it}$ for $\ta_t^\bG$ is defined by 
$P_\bG^{it} (x\hat{1}_h)=\ta_t^\bG(x)\hat{1}_h$ 
for $t\in\R$ and $x\in\lG$. 
Then the following equalities are directly deduced 
from (\ref{eq: W}), (\ref{eq: V}), (\ref{eq: de1}) and (\ref{eq: de2}). 

\begin{align}
W_\bG
=&(\hJ_\bG\oti J_\bG)W_\bG^*(\hJ_\bG\oti J_\bG), 
&
V_\bG
=&(J_\bG\oti\hJ_\bG)V_\bG^*(J_\bG\oti\hJ_\bG),
\label{eq: JW}\\
W_\bG
=&
(P_\bG^{it}\oti P_\bG^{it})W_\bG(P_\bG^{-it}\oti P_\bG^{-it}),
&
V_\bG
=&
(P_\bG^{it}\oti P_\bG^{it})V_\bG(P_\bG^{-it}\oti P_\bG^{-it}),
\label{eq: PW}\\
W_\bG
=&
(P_\bG^{it}\oti \vDe_{h_\bG}^{it})W_\bG(P_\bG^{-it}\oti \vDe_{h_\bG}^{-it}),
&
V_\bG
=&
(\vDe_{h_\bG}^{it}\oti P_\bG^{-it})V_\bG(\vDe_{h_\bG}^{-it}\oti P_\bG^{it}),
\label{eq: vDW}
\\
W_\bG
=&(\hJ_\bG\oti \hJ_\bG)(V_\bG)_{21}^*(\hJ_\bG\oti \hJ_\bG).
&
&
\label{eq: VW}
\end{align}

We denote by $\lG_*^{\rm{fin}}$ 
the set of $\om\in\lG_*$ which satisfies 
$\om(v_{\pi_{i,j}})=0$, $i,j\in I_\pi$ 
all but finite elements $\pi\in \IG$. 
We simply write symbols by omitting $\bG$, 
if no confusion arise. 
For example, we write $\de$ for $\de_\bG$.

\subsection{Discrete quantum groups.}
The notions of a (dual) discrete quantum group 
have been studied in many papers, 
for example, \cite{EfRu}, \cite{KV}, \cite{VanDaele}, \cite{Wor2} 
and \cite{Ya}. 
They have described essentially the same object. 
In this paper, we use a von Neumann algebraic quantum group 
presented in \cite{KV}. 

\begin{defn}
A \emph{discrete quantum group} is a quintuplet 
$(M,\De,\vph,\psi,\vep)$ 
which satisfies the following conditions: 
\begin{enumerate}

\item $M$ is a separable von Neumann algebra. 

\item $\De\col M\ra M\oti M$ is a coproduct, that is, 
it is a faithful normal unital $*$-homomorphism 
satisfying the coassociativity condition, 
\[
(\De\oti\id)\circ\De=(\id\oti\De)\circ\De. 
\]

\item 
$\vph$ is a faithful normal semifinite weight on $M$ 
satisfying the left invariance, 
\[
\vph((\om\oti\id)(\De(x)))=\om(1)\vph(x)\quad 
\mbox{for all}\ \om\in M_*^+,\ x\in m_\vph^+. 
\]

\item 
$\ps$ is a faithful normal semifinite weight on $M$ 
satisfying the right invariance,
\[
\ps((\id\oti\om)(\De(x)))=\om(1)\ps(x)\quad 
\mbox{for all}\ \om\in M_*^+,\ x\in m_\ps^+. 
\]

\item 
$\vep$ is a normal counit, that is, 
it is a normal character on $M$ satisfying 
\[(\vep\oti\id)\circ\De=\id=(\id\oti\vep)\circ\De.
\] 

\end{enumerate}
\end{defn}

In fact, those weights $\vph$, $\psi$ are uniquely determined 
up to scalar multiplications and the counit $\vep$ is unique. 
From now, we simply write $(M,\De)$ for $(M,\De,\vph,\psi,\vep)$ 
by omitting $\vph$, $\ps$ and $\vep$ 
once they are given. 

For a compact quantum group $\bG$, we construct the dual discrete 
quantum group as follows (see for example \cite{BS}, \cite{KV}). 
We define the \emph{left group algebra} 
and \emph{right group algebra} by
\[
\LG=\ovl{\{(\om\oti\id)(W_\bG)\mid \om\in \lG_*\}}^{\,\rm{w}}, 
\ 
\RG=\ovl{\{(\id\oti\om)(V_\bG)\mid \om\in \lG_*\}}^{\,\rm{w}}.
\] 
By (\ref{eq: decomp}), 
$\LG$ and $\RG$ are generated by $\{f_{\pi_{i,j}}\}$ 
and $\{e_{\pi_{i,j}}\}$, respectively. 
Hence they are isomorphic to 
the von Neumann algebra direct sum of 
the matrix algebras $\{B(H_\pi)\}_{\pi\in\IG}$. 
Set 
\[
\LG^{\rm{fin}}=\{(\om\oti\id)(W_\bG)\mid \om\in \lG_*^{\rm{fin}}\},
\ 
\RG^{\rm{fin}}=\{(\id\oti\om)(V_\bG)\mid \om\in \lG_*^{\rm{fin}}\}. 
\]
They are algebraic direct sum of the matrix algebras. 
We note the commutant property $\LG'=\RG=\hJ \LG \hJ$. 
We define a minimal central projection corresponding to 
$\pi\in\IG$ by $1_\pi=\sum_{i\in I_\pi}f_{\pi_{i,i}}$. 
In particular, $e_\btr=1_\btr$. 
We denote by $\LG_*^{\rm{fin}}$ the subset of $\LG_*$ 
which consists of $\om\in \LG_*$ such that 
$\om(f_{\pi_{i,j}})=0$, $i,j\in I_\pi$ all but finite elements $\pi\in\IG$. 

We prepare the coproducts $\De_L$ and $\De_R$ defined by 
\[
\De_L(x)=W_\bG(x\oti 1)W_\bG^*
\quad \mbox{for}\ 
x\in \LG,
\]
\[
\hspace{1pt}
\De_R(x)=V_\bG^*(1\oti x)V_\bG
\hspace{19pt} \mbox{for}\ 
x\in \RG. 
\]
We define the normal counit $\hvep\col \LG\ra \C$ by 
\[
\hvep(x)e_\btr=x e_\btr\quad
\mbox{for}\ x\in\LG. 
\]
In fact, there exist left, right invariant weights on $\LG$ and $\RG$ 
\cite{KV}, and they are discrete quantum groups. 
Note that every discrete quantum group arises 
as the left (right) group algebra of a compact quantum group 
(see \cite{BS, KV} for duality theory). 
In this paper, we use the symbol $\bhG$ for the discrete quantum group 
$(\RG,\De_R)$. 
We simply write $\De_\bG$ (or $\De$) 
for $\De_L$ and $\De_R$ when it is not ambiguous. 

We define a positive operator $F$ affiliated with $\LG$ by 
\[
F=\sum_{\pi\in \IG}\sum_{i\in I_\pi}F_{\pi_{i,i}}f_{\pi_{i,i}}. 
\]
Using (\ref{eq: fpi}), for all $t\in\R$ and $x\in\AG$ we have 
\[
F^{it} (x \hat{1}_h)=(f_{-it}\oti\id)(\de(x))\hat{1}_h,
\quad
\hJ F^{it} \hJ (x \hat{1}_h)= (\id\oti f_{it})(\de(x))\hat{1}_h. 
\]
Then we have 
\begin{equation}\label{eq: vDP}
\vDe_h^{it}=F^{-it} \hJ F^{it} \hJ,
\quad
P^{it}=F^{-it} \hJ F^{-it} \hJ. 
\end{equation}

The antipode $\hS$ on $\LG$ is defined as follows. 
Since the map $\lG_*\ni\om \mapsto (\om\oti\id)(W_\bG)\in\lG$ is 
injective, we can define an invertible antimultiplicative map 
$\hS\col \LG^{\rm{fin}}\ra \LG^{\rm{fin}}$ by 
\[
\hS((\om\oti\id)(W_\bG))=(\om\oti\id)(W_\bG^*)
\quad\mbox{for all}
\ \om\in\lG_*^{\rm{fin}}. 
\]
We define the unitary antipode $\hR$ and 
the scaling automorphism group $\{\hta_t\}_{t\in\R}$ 
on $\LG$ by 
\begin{equation}\label{eq: hRhta}
\hR(x)=J x^* J,\quad 
\hta_t(x)=F^{it}x F^{-it}
\quad\mbox{for}\ 
x\in\LG,\ t\in\R. 
\end{equation}
Using (\ref{eq: JW}), (\ref{eq: PW}) and (\ref{eq: vDP}), 
we have $\hS=\hR\circ \hta_{-i/2}$ on $\LG^{\rm{fin}}$. 
The relations among $\De$, $\hR$ and $\hta$ are as follows. 
\begin{equation}\label{eq: DehR}
\De\circ \hR=(\hR\oti\hR)\circ\De^\op,
\quad
\De\circ\hta_t=(\hta_t\oti\hta_t)\circ\De. 
\end{equation}

Take $\opi\in\IG$ satisfying $\hS(1_\pi)=1_\opi$. 
Let 
$\mu\col \LG^{\rm{fin}}\oti \LG^{\rm{fin}}\ra \LG^{\rm{fin}}$ 
be the multiplication map. 
We prove the well-known results for readers' convenience. 

\begin{lem}\label{lem: anti-couni}
For any $x\in\LG$ and $\pi\in\IG$, one has
\[
\mu\big{(}(\id\oti\hS)({}_\opi\De_\pi(x))\big{)}
=\hvep(x)1_\opi
=\mu\big{(}(\hS\oti\id)({}_\pi\De_\opi(x))\big{)}, 
\]
where ${}_\pi\De_\opi(x)=\De(x)(1_\pi\oti1_\opi)$.
\end{lem}
\begin{proof}
We prove 
$\mu\big{(}(\id\oti\hS)({}_\opi\De_\pi(x))\big{)}
=\hvep(x)1_\opi$.  
Since 
${}_\pi\De_\opi(1_\rho)$ 
is non-zero for finite elements $\rho\in\IG$, 
we may assume that $x=(\om\oti\id)(W_\bG)$, 
$\om\in\lG_*^{\rm{fin}}$. 
Then we have 
\begin{align*}
(\id \oti\hS)({}_\opi\De_\pi(x))
=&\,
(\id \oti\hS)
({}_\opi\De_\pi((\om\oti\id)(W_\bG)))
\\
=&\,
(\id \oti\hS)
((\om\oti\id\oti\id)((\id\oti{}_\opi\De_\pi)(W_\bG)))
\\
=&\,
(\id \oti\hS)
((\om\oti\id\oti\id)((W_\bG)_{12}(W_\bG)_{13}(1\oti1_\opi\oti1_\pi)))
\\
=&\,
(\om\oti\id\oti\id)((W_\bG)_{12}(W_\bG)_{13}^*(1\oti1_\opi\oti1_\opi)).
\end{align*}
Hence 
\begin{align*}
\mu\big{(}(\id \oti\hS)({}_\opi\De_\pi(x))\big{)}
=&\,
\mu\big{(}
(\om\oti\id\oti\id)((W_\bG)_{12}(W_\bG)_{13}^*(1\oti1_\opi\oti1_\opi))
\big{)}
\\
=&\,
(\om\oti\id)(W_\bG W_\bG^*(1\oti 1_\opi))
\\
=&\,
\om(1)1_\opi=\hvep(x)1_\opi. 
\end{align*}
Similarly we can prove 
$\mu\big{(}(\hS\oti\id)({}_\pi\De_\opi(x))\big{)}=\hvep(x)1_\opi$. 
\end{proof}

\begin{lem}\label{lem: anti-repn}
Let $K$ be a Hilbert space and $\mathscr{V}\in B(K)\oti \LG$ 
a unitary representation of $(\LG,\De)$, that is, 
$\meV$ is a unitary satisfying 
$(\id\oti\De)(\meV)=\meV_{12}\meV_{13}$. 
Then for any $\pi\in\IG$, one has 
\[
(\id\oti \hS)(\meV_\pi)=\meV_\opi^*
\]
where $\meV_\pi=\meV(1\oti 1_\pi)$.  
\end{lem}
\begin{proof}
Set $x=(\om\oti\id)(\meV)\in\LG$, $\om\in B(K)_*$. 
For any $\pi\in\IG$, we have 
\begin{align*}
(\id\oti\hS)({}_\opi \De_{\pi}(x))
=&\,
(\id\oti\hS)\big{(}
(\om\oti\id\oti\id)((\id\oti{}_\opi\De_\pi)(\meV))\big{)}
\\
=&\,
(\om\oti\id\oti\hS)((\meV_{\opi})_{12}(\meV_{\pi})_{13})
\\
=&\,
(\om\oti\id\oti\id)
\big{(}(\meV_{\opi})_{12}(\id\oti\hS)(\meV_{\pi})_{13}\big{)}
\end{align*}

By the previous lemma, we have 
$\mu\big{(}(\id\oti \hS)({}_\opi\De_\pi(x))\big{)}=\hvep(x)1$. 
This implies 
$(\om\oti\id)(\meV_\opi(\id\oti\hS)(\meV_{\pi}))=\om(1)$ for any 
$\om\in B(K)_*$. 
Hence $(\id\oti\hS)(\meV_{\pi})=\meV_\opi^*$. 
\end{proof}

\subsection{Amenability.}
We recall the notion of amenability of a discrete quantum group. 
For a detail of the theory, 
readers are referred to 
\cite{BeCoTu}, \cite{BeMuTu1}, \cite{BeMuTu2}, \cite{To1} 
and references therein. 
A discrete quantum group $(M,\De)$ is said to be \textit{amenable} 
if there exists a state $m\in M^*$ such that 
$m((\om\oti\id)(\De(x)))=\om(1)m(x)$ 
for all $\om\in M_*$ and $x\in M$. 
The state $m$ is called a \emph{left invariant mean}. 
It is known that $(\LG,\De_L)$ is amenable if and only if 
$(\RG,\De_R)$ is amenable. 
If it is the case, 
we say that $\bG$ is \textit{co-amenable}. 
Note that the counit $\vep_\bG$ is norm-bounded 
if and only if $\bG$ is co-amenable 
\cite[Theorem 4.7]{BeMuTu2}, \cite[Theorem 3.8]{To1}. 
The amenability is also equivalent to the universality of $\CG$, that is, 
for any $C^*$-algebra $B$, 
any $*$-homomorphism $\si\col \AG\ra B$ 
extends to a $*$-homomorphism $\si\col \CG \ra B$ 
\cite[Theorem 3.6]{BeMuTu1}. 

\subsection{Right $\bG$-action $\al$ and left $\bhG$-action $\be$.}
\label{subsect: al-be}
We prepare two maps $\al\col B(\ltG)\ra B(\ltG)\oti\lG$ 
and $\be\col B(\ltG)\ra \RG\oti B(\ltG)$ 
which will be frequently used in our study. 
A map $\al$ is defined by 
\[
\al(x)=V_\bG (x\oti 1) V_\bG^*
\quad 
\mbox{for}\ x\in B(\ltG). 
\]
Then $\al$ is a right action of $\bG$ on $B(\ltG)$, 
that is, $(\al\oti\id)\circ\al=(\id\oti\de_\bG)\circ\al$ holds. 
The conditional expectation $E_\al=(\id\oti h)\circ\al$ 
maps $B(\ltG)$ onto 
$\LG=\{x\in B(\ltG)\mid \al(x)=x\oti1\}$. 
The other one, $\be$ is defined by 
\[
\be(x)=V_\bG^*(1\oti x)V_\bG 
\quad 
\mbox{for}\ x\in B(\ltG). 
\]
Then $\be$ is a left action of $\bhG$ on $B(\ltG)$, 
that is, $(\id\oti\be)\circ\be=(\De_R\oti\id)\circ\be$ holds. 
We also call $\be$ a \emph{coaction} of $\bG$. 
Note that both the actions $\al$, $\be$ preserve 
not only $\lG$ but also $\RG$. 

\subsection{Right coideals.}
We introduce the notion of a right coideal. 
Our basic references for theory of right coideals 
are \cite{ILP} and \cite{To2}. 

\begin{defn}
Let $\bG$ be a compact quantum group. 
Let $B\subs \lG$ be a von Neumann subalgebra. 
We say that 
$B$ is a \textit{right coideal} 
if 
$\de(B)\subs B\oti \lG$. 
\end{defn}

A left coideal in $\LG$ is similarly defined. 
For a right coideal $B\subs \lG$ and a left coideal $C\subs \LG$, 
we define 
\[
\hB=B'\cap \LG,\quad \hC=C'\cap \lG. 
\]
We prove a quantum group 
version of \cite[Theorem 4.6]{ILP} as follows. 

\begin{lem}\label{lem: lattice}
Let $\bG$ be a compact quantum group. 
Let $B\subs \lG$ be a right coideal and $C\subs\LG$ a left coideal. 

\begin{enumerate}

\item 
$\hB\subs \LG$ is a left coideal and $\hC\subs \lG$ is a right coideal.

\item 
The map $B\mapsto\hB$ is a lattice anti-isomorphism between the set of 
right coideal of $\lG$ and left coideal of $\LG$. 
The inverse map is given by $C\mapsto \hC$. 
\end{enumerate}
\end{lem}
\begin{proof}
(1). 
It is similarly proved as in \cite[Theorem 4.6]{ILP}. 

(2). 
First we show that the map $C\mapsto \hC$ is injective. 
Let $C\subs\LG$ be a left coideal. 
Then $\hC'=(C'\cap\lG)'=C\vee \lG'$. 
We can adapt the proof of \cite[Theorem 4.6 (ii)]{ILP} 
to deduce $\hC'=\ovl{C\, \lG'}^{\,\rm{w}}$. 
Recall the right $\bG$-action $\al$ on $B(\ltG)$ 
and the conditional expectation $E_\al\col B(\ltG)\ra\LG$. 
We apply $E_\al$ to $J\hC'J=\ovl{JCJ\,\lG}^{\,\rm{w}}$. 
Since $JCJ\subs J\LG J\subs\LG$ and 
the restriction of $E_{\al}$ to $\lG$ is the Haar state $h$, 
we have $E_{\al}(J\hC'J)=JCJ$. 
Hence the map $C\mapsto \hC$ is injective. 
Set $B=\hC$. 
Then $\hB'\cap\lG=\hC$ holds. 
The injectivity of the map $C\mapsto \hC$ yields $\hB=C$. 

Next we show that the map $B\mapsto \hB$ is injective. 
Let $B_i\subs\lG$, $i=1,2$ be right coideals. 
Then $\hB_i'=(B_i'\cap \LG)'=B_i\vee \RG$. 
Let $B_i\rti_\de \bG=\de(B_i)\vee \C\oti\RG$ be the crossed product. 
Then we have 
$W_\bG (B_i\rti_\de \bG) W_\bG^*=\C\oti (B_i\vee \RG)=\C\oti \hB_i'$. 
Hence $\hB_1=\hB_2$ if and only if 
$B_1\rti_\de \bG=B_2\rti_\de \bG$. 
Then the biduality theorem 
(see \cite[Th\'{e}or\`{e}me 7.5]{BS} or \cite[Theorem 2.6]{Vaes1}) implies 
that it is equivalent to 
$B_1\oti B(\ltG)=B_2\oti B(\ltG)$, and $B_1=B_2$. 
Hence the map $B\mapsto \hB$ is injective. 
\end{proof}

\subsection{Quantum subgroups.}
For the definition of quantum subgroups, 
we follow that of \cite{Po2}, 
in which matrix pseudogroups \cite{Wor1} are treated. 
\begin{defn}
Let $\bG$ and $\bH$ be compact quantum groups. 
\begin{enumerate}

\item 
Suppose that there exists a surjective $*$-homomorphism 
$r_\bH\col \AG\ra \AH$ such that 
$\de_\bH\circ r_\bH=(r_\bH\oti r_\bH)\circ\de_\bG$. 
Then we say that the pair $\{\bH,r_\bH\}$ 
is an \textit{algebraic quantum subgroup} of $\bG$. 

\item
Suppose that there exists a surjective $*$-homomorphism 
$r_\bH\col C(\bG)\ra C(\bH)$ such that 
$\de_\bH\circ r_\bH=(r_\bH\oti r_\bH)\circ\de_\bG$. 
Then we say that the pair $\{\bH,r_\bH\}$ 
is a \textit{quantum subgroup} of $\bG$. 

\end{enumerate}
\end{defn}
In the above cases, the map $r_\bH$ is called a \emph{restriction map}. 
We present basic properties on (algebraic) quantum subgroups 
in the following lemmas. 
Those are probably well-known for specialists, 
but we prove them for the sake of readers' convenience. 

\begin{lem}\label{lem: algsub}
Let $\bG$ and $\bH$ be compact quantum groups. 
\begin{enumerate}

\item
If $\{\bH,r_\bH\}$ is a quantum subgroup of $\bG$, 
then it is an algebraic quantum subgroup of $\bG$. 

\item
Suppose $\bG$ is co-amenable. 
Then any algebraic quantum subgroup of $\bG$ is 
naturally regarded as a quantum subgroup. 

\end{enumerate}
\end{lem}
\begin{proof}
(1). 
Let $\pi\in\IG$. 
Then $(\id\oti r_\bH)(v_\pi)\in B(H_\pi)\oti \CH$ 
is a finite dimensional unitary representation of $\bH$. 
Hence 
it is a finite direct sum of irreducible representations of $\bH$. 
This implies $(\id\oti r_\bH)(v_\pi)\in B(H_\pi)\oti \AH$. 
Hence $r_\bH(\AG)\subs\AH$. 
Next we show the converse inclusion. 
Take any $\rho\in\IH$. 
Let $w_\rho\in B(K_\rho)\oti \AH$ be a corresponding 
irreducible unitary representation. 
Consider the bounded linear map $\th\col \CG \ra B(K_\rho)$ 
defined by 
$\th(x)=(\id\oti h_\bH)(w_\rho^* (1\oti r_\bH(x)))$. 
Since $r_\bH$ is surjective, $\th$ is a non-zero map. 
By density of $\AG\subs\CG$, 
There exists $\pi\in \IG$ such that 
$(\id\oti\th)(v_\pi)\neq 0$. 
This shows that the unitary representation $(\id\oti r_\bH)(v_\pi)$ 
contains $w_\rho$. 
Hence all the entries of $w_\rho$ are contained in $r_\bH(\AG)$, 
and $\AH\subs r_\bH(\AG)$. 

(2). 
Let $r_\bH\col \AG\ra \AH$ be a restriction map. 
Since $\CG$ is universal, 
the map extends to $r_\bH\col \CG\ra\CH$. 
The image contains a total subspace $\AH$ in $\CH$, 
and $r_\bH$ is surjective. 
By continuity of $r_\bH$, 
the relation $\de_\bH\circ r_\bH=(r_\bH\oti r_\bH)\circ\de_\bG$ 
holds on $\CG$. 
Hence $\{\bH,r_\bH\}$ is a quantum subgroup of $\bG$. 
\end{proof}

\begin{lem}\label{lem: algsub-scaling}
Let $\{\bH,r_\bH\}$ be an algebraic quantum subgroup of $\bG$. 
On $\AG$,  
\begin{enumerate}

\item $\vep_\bH\circ r_\bH=\vep_\bG$, 

\item $\ta_t^\bH\circ r_\bH=r_\bH\circ \ta_t^\bG$ for all $t\in\R$, 

\item $R_\bH\circ r_\bH= r_\bH\circ R_\bG$. 

\end{enumerate}
\end{lem}
\begin{proof}
(1). 
Take any $\pi\in\IG$. 
Since $(\id\oti r_\bH)(v_\pi)$ is a unitary representation of $\bH$, 
$(\id\oti \vep_\bH\circ r_\bH)(v_\pi)=1=(\id\oti\vep_\bG)(v_\pi)$. 
This implies $\vep_\bH\circ r_\bH=\vep_\bG$. 

(2). 
Take any $\pi\in\IG$. 
Since $(\id\oti r_\bH)(v_\pi)$ is a unitary representation of $\bH$, 
$(\id\oti \ka_\bH\circ r_\bH)(v_\pi)=(\id\oti r_\bH)(v_\pi)^*
=(\id\oti r_\bH\circ \ka_\bG)(v_\pi)$. 
This implies $r_\bH\circ\ka_\bG=\ka_\bH\circ r_\bH$. 
Since $\ka_\bG^2=\ta_{-i}^\bG$ and $\ka_\bH^2=\ta_{-i}^\bH$, 
we have $r_\bH\circ \ta_{-i}^\bG=\ta_{-i}^\bH\circ r_\bH$. 
Let $\pi\in\IG$ and set $w=(\id\oti r_\bH)(v_\pi)$. 
Then we have 
\begin{align*}
(F_\pi^\bG\oti1)w((F_\pi^\bG)^{-1}\oti1)
\!=&
(\id\oti r_\bH)
\big{(}(F_\pi^\bG\oti1)v_\pi((F_\pi^\bG)^{-1}\oti1)\big{)}
\!=\!
(\id\oti r_\bH\circ \ta_{-i}^\bG)(v_\pi)
\\
=&
(\id\oti \ta_{-i}^\bH\circ r_\bH)(v_\pi)
=
(\id\oti\ta_{-i}^\bH)(w)
\\
=&
(F_w^\bH\oti 1)w((F_w^\bH)^{-1}\oti1). 
\end{align*}
This shows $k_\pi:=(F_w^\bH)^{-1} F_\pi^\bG\in B(H_\pi)$ 
is an intertwiner of $w$, 
and in particular, $F_w^\bH k_\pi=k_\pi F_w^\bH$ holds. 
Hence $F_\pi^\bG=F_w^\bH k_\pi=k_\pi F_w^\bH$. 
Since the positive operators $F_\pi^\bG$ and $(F_w^\bH)^{-1}$ 
commute, 
$k_\pi$ is positive. 
In particular, we have $F_\pi^{it}=(F_w^{\bH})^{it} k_\pi^{it}$ 
for all $t\in\R$. 
Hence for any $\pi\in\IG$, we have
\begin{align*}
(\id\oti r_\bH\circ \ta_t^\bG)(v_\pi)
=&\,
(\id\oti r_\bH)
\big{(}((F_\pi^\bG)^{it}\oti 1) v_\pi ((F_\pi^\bG)^{-it}\oti 1)\big{)}
\\
=&\,
((F_\pi^\bG)^{it}\oti 1) w ((F_\pi^\bG)^{-it}\oti 1)
\\
=&\,
((F_w^\bH)^{it} k_\pi^{it} \oti 1)w
(k_\pi^{-it}(F_w^\bH)^{-it}\oti 1)
\\
=&\,
((F_w^\bH)^{it} \oti 1)w
((F_w^\bH)^{-it}\oti 1)
\\
=&\,
(\id\oti \ta_t^\bH)(w)
\\
=&\,
(\id\oti \ta_t^\bH\circ r_\bH)(v_\pi). 
\end{align*}
Therefore the desired relation holds. 

(3). 
It follows from $r_\bH\circ\ka_\bG=\ka_\bH\circ r_\bH$ and (2). 
\end{proof}

Let $\{\bH,r_\bH\}$ be 
an algebraic quantum subgroup of $\bG$. 
Since $W_\bG(1\oti 1_\pi)\in \AG\oti \LG 1_\pi$ for any $\pi\in\IG$, 
we can define a left unitary representation 
of $\bH$, $(r_\bH\oti\id)(W_\bG)\in \lH\oti \LG$ by 
\[
(r_\bH\oti\id)(W_\bG)=\sum_{\pi\in\IG}(r_\bH\oti\id)(W_\bG(1\oti1_\pi)). 
\]

\begin{lem}\label{lem: algsub-corresp}
Let $\bG$ and $\bH$ be compact quantum groups. 
\begin{enumerate}

\item 
Assume that $\{\bH,r_\bH\}$ is 
an algebraic quantum subgroup of $\bG$. 
Then there exists a faithful normal unital $*$-homomorphism 
$r_\bH^t\col \LH\ra \LG$ such that 
\[
(\id\oti r_\bH^t)(W_\bH)=(r_\bH\oti\id)(W_\bG). 
\]

\item 
Assume that there exists a faithful normal unital $*$-homomorphism 
$\iota\col \LH\ra \LG$ such that 
\[
\De_\bG\circ\iota= (\iota\oti\iota)\circ \De_\bH. 
\]
Then there exists a restriction map $r_\bH\col \AG\ra\AH$ 
such that $\iota=r_\bH^t$. 

\end{enumerate}
\end{lem}
\begin{proof}
(1). 
Since the map $\lH_* \ni \om\mapsto (\om\oti\id)(W_\bH)\in \LH$ 
is injective, 
we can define a linear map $r_\bH^t\col \LH^{\rm{fin}}\ra \LG$ by 
\[
r_\bH^t((\om\oti\id)(W_\bH))=(\om\oti\id)((r_\bH\oti\id)(W_\bG))
\quad\mbox{for}\ 
\om\in \lH_{*}^{\rm{fin}}. 
\]
It is easy to verify the multiplicativity of $r_\bH^t$ 
by the pentagon equalities (\ref{eq: Wpenta}) on $W_\bH$ and $W_\bG$. 
We show that $r_\bH^t$ preserves the involutions. 
Take any $\om\in\lH_*^{\rm{fin}}$ and $\pi\in\IG$. 
Using the equality $r_\bH\circ\ka_\bG=\ka_\bH\circ r_\bH$, we have 
\begin{align*}
r_\bH^t((\om\oti\id)(W_\bH))^*(1\oti1_\pi)
=&\,
\big{(}(\om\oti\id)( (r_\bH\oti\id)(W_\bG))\big{)}^*(1\oti1_\pi)
\\
=&\,
(\ovl{\om}\circ r_\bH\oti\id)(W_\bG^*(1\oti1_\pi))
\\
=&\,
(\ovl{\om}\circ r_\bH\oti\id)((\ka_\bG\oti\id)(W_\bG)(1\oti1_\pi))
\\
=&\,
(\ovl{\om}\circ\ka_\bH\oti\id)\big{(}(r_\bH\oti\id)(W_\bG)(1\oti1_\pi)\big{)}
\\
=&\,
r_\bH^t((\ovl{\om}\circ\ka_\bH\oti\id)(W_\bH))(1\oti 1_\pi)
\\
=&\,
r_\bH^t\big{(}\big{(}(\om\oti\id)(W_\bH)\big{)}^*\big{)}(1\oti 1_\pi). 
\end{align*}
Hence $r_\bH^t$ is a $*$-homomorphism. 
Since $\LH$ is a von Neumann algebra direct sum of 
$\{B(H_\rho)\}_{\rho\in\IH}$, 
the map $r_\bH^t$ extends to a normal $*$-homomorphism from $\LH$ to $\LG$. 
Then the desired equality holds. 
We show that $r_\bH^t$ is unital. 
Take any $\pi\in\IG$ and $\rho\in\IH$. 
Let $N_{\pi|_\bH}^\rho$ be the multiplicity of $\rho$ in the unitary 
representation $(r_\bH\oti\id)(W_\bG(1\oti1_\pi))$. 
Then the map $r_\bH^t\col \LH1_\rho\ra \LG1_\pi$ is 
the $N_{\pi|_\bH}^\rho$ times amplification, and $r_\bH^t$ is unital. 
Since $r_\bH$ is surjective, $r_\bH^t$ is injective. 

(2). 
Since the linear map $\LG_*\ni \om\mapsto (\id\oti\om)(W_\bG)\in \lG$ 
is injective, we can define a linear map $r_\bH\col \AG\ra \AH$ by 
\[
r_\bH((\id\oti\om)(W_\bG))=(\id\oti \om\circ \iota)(W_\bH)
\quad\mbox{for}\ \om\in \LG_*^{\rm{fin}}.
\]
First we show that $r_\bH$ is multiplicative. 
Take any $\om,\th\in \LG_*^{\rm{fin}}$. 
Then $(\om\oti\th)\circ\De_\bG\in \LG_*^{\rm{fin}}$, 
and 
\begin{align*}
r_\bH((\id\oti\om)(W_\bG))\,r_\bH((\id\oti\th)(W_\bG))
=&\,
(\id\oti \om\circ \iota)(W_\bH)(\id\oti \th\circ \iota)(W_\bH)
\\
=&\,
(\id\oti \om\circ \iota\oti \th\circ\iota)((W_\bH)_{12}(W_\bH)_{13})
\\
=&\,
(\id\oti \om\circ \iota\oti \th\circ\iota)((\id\oti\De_\bH)(W_\bH))
\\
=&\,
(\id\oti\om\oti\th)\big{(}(\id\oti \De_\bG)((\id\oti\iota)(W_\bH))\big{)}
\\
=&\,
r_\bH((\id\oti (\om\oti\th)\circ\De_\bG)(W_\bG))
\\
=&\,
r_\bH((\id\oti\om)(W_\bG)(\id\oti\th)(W_\bG)). 
\end{align*}
Next we show $r_\bH$ preserves the involutions. 
Let $\pi\in\IG$. 
Take $\om\in \LG_*^{\rm{fin}}$ such that $\om$ is 
equal to zero on $B(H_\rho)$ if $\rho\neq \opi$. 
Since $\De_\bG\circ\iota=(\iota\oti\iota)\circ\De$, 
$(\id\oti\iota)(W_\bH)$ is a unitary representation of 
the discrete quantum group $(\LG,\De)$. 
Set $\meV=(\id\oti\iota)(W_\bH)$. 
By Lemma \ref{lem: anti-repn}, we have 
$(\id\oti\hS)(\meV_\pi)=\meV_\opi^*$. 
Then 
\begin{align*}
r_\bH((\id\oti\om)(W_\bG))^*
=&\,
\big{(}(\id\oti \om\circ \iota)(W_\bH)\big{)}^*
=
(\id\oti \ovl{\om}\circ \iota)(W_\bH^*)
\\
=&\,
(\id\oti\ovl{\om})(\meV_\opi^*)
=
(\id\oti\ovl{\om}\circ\hS)(\meV_\pi)
\\
=&\,
r_\bH((\id\oti \ovl{\om}\circ\hS)((W_\bG)_\pi))
=
r_\bH\big{(}(\id\oti \ovl{\om})((W_\bG)_\opi^*)\big{)}
\\
=&\,
r_\bH\big{(}((\id\oti\om)(W_\bG))^*\big{)}. 
\end{align*}
By taking summations on $\om$, we see that the above equality holds 
for any $\om\in\LG_*^{\rm{fin}}$. 
Hence $r_\bH$ is a unital $*$-homomorphism. 
Let $\mu\in \LH_*^{\rm{fin}}$. 
Then there exists $\tilde{\mu}\in \LG_*^{\rm{fin}}$ such that 
$\mu=\tilde{\mu}\circ \iota$. 
Then we have 
$r_\bH((\id\oti\tilde{\mu})(W_\bG))=(\id\oti \mu)(W_\bH)$. 
Hence $r_\bH$ is surjective. 
Finally we show that $r_\bH$ is a restriction map. 
Let $\om\in\LG_*^{\rm{fin}}$. 
Then 
\begin{align*}
\de_\bH\big{(}r_\bH((\id\oti\om)(W_\bG))\big{)}
=&\,
\de_\bH((\id\oti\om\circ \iota)(W_\bH))
\\
=&\,
(\id\oti\id\oti \om\circ \iota)((W_\bH)_{13}(W_\bH)_{23})
\\
=&\,
(r_\bH\oti r_\bH)
\big{(} 
(\id\oti\id\oti \om\circ \iota)((W_\bG)_{13}(W_\bG)_{23})
\big{)}
\\
=&\,
(r_\bH\oti r_\bH)\circ\de_\bG
((\id\oti\om)(W_\bG)). 
\end{align*}
Hence $\{\bH,r_\bH\}$ is an algebraic quantum subgroup of $\bG$. 
It is clear that $\iota=r_\bH^t$ by definition of $r_\bH$. 
\end{proof}

On heredity of co-amenability to quantum subgroups, 
we have the following lemma. 

\begin{lem}\label{lem: heredity}
Let $\bG$ be a compact quantum group. 
Then the following statements are equivalent: 
\begin{enumerate}

\item 
$\bG$ is co-amenable. 

\item 
Any quantum subgroup of $\bG$ is co-amenable. 

\item 
There exists a co-amenable quantum subgroup of $\bG$. 

\end{enumerate}
\end{lem}
\begin{proof}
The implication (2) to (3) is trivial. 
Suppose (3) holds. 
Let $\{\bH,r_\bH\}$ be a co-amenable quantum subgroup. 
Then $\vep_\bG=\vep_\bH\circ r_\bH$ is bounded, and (1) holds. 
We have to show that (1) implies (2). 
Let $\{\bH,r_\bH\}$ be a quantum subgroup of $\bG$. 
By the previous lemma, there exists an embedding map 
$r_\bH^t\col \LH\ra \LG$ preserving the coproducts. 
Let $m\in \LG^*$ be a left invariant mean. 
Then the state $m\circ r_\bH^t \in \LH^*$ 
is also a left invariant mean on $\LH$. 
Hence the discrete quantum group $(\LH,\De)$ is amenable, 
and $\bH$ is co-amenable. 
\end{proof}

Let $\{\bH,r_\bH\}$ be a quantum subgroup of $\bG$. 
We define a map $\ga_\bH\col \CG\ra \CH\oti\CG$ by 
$\ga_\bH(x)=(r_\bH\oti\id)(\de_\bG(x))$ for $x\in\CG$. 
Then $\ga_\bH$ is a left action of $\bH$ on $\CG$, that is, 
$\ga_\bH$ satisfies 
$(\id\oti\ga_\bH)\circ\ga_\bH=(\de_\bH\oti\id)\circ\ga_\bH$. 
We set the fixed point algebra 
$C(\bH\setm\bG)=\{x\in\CG\mid \ga_\bH(x)=1\oti x\}$. 
Since the action $\ga_\bH$ preserves the Haar state of $\bG$, 
it extends to an action on $\lG$. 
We denote the fixed point algebra by $L^\infty(\bH\setm\bG)$ 
as well. 
Then $L^\infty(\bH\setm\bG)$ is the range of 
the conditional expectation 
$E_\bH=(h_\bH\oti\id)\circ\ga_\bH$. 

We often identify $\LH$ with the subalgebra of $\LG$ 
via the inclusion map $r_\bH^t\col \LH\ra \LG$ defined 
in Lemma \ref{lem: algsub-corresp}. 
By the identification, we have $W_\bH=(r_\bH\oti\id)(W_\bG)$. 
Similarly, we can do 
$\RH\subs\RG$ and $(\id\oti r_\bH)(V_\bG)=V_\bH$. 
Then the left action of $\bH$ on $\lG$ is given by 
$\ga_\bH(x)=W_\bH^*(1\oti x)W_\bH$ for $x\in\lG$. 

\begin{lem}\label{lem: above}
In the above setting, one has 
\[
\LH'\cap \lG=L^\infty(\bH\setm\bG), 
\quad
L^\infty(\bH\setm\bG)'\cap \LG=\LH. 
\] 
\end{lem}
\begin{proof}
By definition of $r_\bH$, 
$x\in \LH'\cap \lG$ if and only if 
$W_\bH^*(1\oti x)W_\bH=1\oti x$. 
It is equivalent to $x\in L^\infty(\bH\setm\bG)$. 
Hence $\LH'\cap \lG=L^\infty(\bH\setm\bG)$. 
By Lemma \ref{lem: lattice}, we obtain the second equality. 
\end{proof}

\begin{lem}\label{lem: subs}
Let $\{\bH,r_\bH\}$ and $\{\bK,r_\bK\}$ 
be quantum subgroups of $\bG$. 
Assume that $L^\infty(\bH\setm \bG)=L^\infty(\bK\setm\bG)$. 
Then there exists a $*$-isomorphism 
$\th \col \CH\ra \CK$ such that 
$r_\bK=\th\circ r_\bH$ on $\CG$. 
\end{lem}
\begin{proof}
By the previous lemma, we have $r_\bH^t(\LH)=r_\bK^t(\LK)$. 
Applying Lemma \ref{lem: algsub-corresp} to left group algebra 
$r_\bH^t(\LH)$, 
we have a $*$-isomorphism $\th\col \CH\ra\CK$ such that 
$(\th\oti\id)((\id\oti r_\bH^t)(W_\bH))
=(\id\oti r_\bK^t)(W_\bK)$. 
This implies 
$(\th\circ r_\bH\oti\id)(W_\bG)=(r_\bK\oti\id)(W_\bK)$, 
and $r_\bK=\th\circ r_\bH$ on $\CG$. 
\end{proof}

\section{Right coideals of quotient type}

\begin{defn}
Let $B\subs \lG$ be a right coideal. 
We say that 
\begin{enumerate}

\item 
$B$ is of \textit{quotient type} 
if there exists a quantum subgroup $\{\bH,r_\bH\}$ of $\bG$ 
such that $B=L^\infty(\bH\setm \bG)$, 

\item 
$B$ 
has the \textit{expectation property} 
if there exists a faithful normal conditional expectation 
$E_B\col \lG\ra B$ such that $h\circ E_B=h$, 

\item 
$B$ 
has the \textit{coaction symmetry} 
if $\be(B)\subs \RG\oti B$. 

\end{enumerate}
\end{defn}

Typical examples of right coideals are 
given by taking quotients. 
In fact, they have the expectation property and the coaction symmetry 
as follows. 

\begin{lem}\label{lem: exco}
A right coideal of quotient type 
has 
the expectation property and the coaction symmetry. 
\end{lem}
\begin{proof}
Let $\{\bH,r_\bH\}$ be a quantum subgroup of $\bG$. 
Set $B=L^\infty(\bH\setm\bG)$. 
It is easy to see that the conditional expectation 
$E_\bH$ preserves the Haar state $h$. 
Next we verify the coaction symmetry. 
Let $x\in C(\bH\setm\bG)$. 
Then we have 
\begin{align*}
(\id\oti\ga_\bH)(\be(x))
=&\,
(\id\oti (r_\bH\oti\id)\circ\de)(\be(x))
\\
=&\,
(\id\oti (r_\bH\oti\id)\circ\de)(V_\bG^*(1\oti x)V_\bG)
\\
=&\,
(\id\oti r_\bH\oti\id)
((V_\bG)_{13}^* (V_\bG)_{12}^* (1\oti\de(x)) 
(V_\bG)_{12}(V_\bG)_{13})
\\
=&\,
(V_\bG)_{13}^* (V_\bH)_{12}^* (1\oti1\oti x)
(V_\bH)_{12}(V_\bG)_{13}
\\
=&\,
(V_\bG)_{13}^* (1\oti1\oti x)(V_\bG)_{13}
\\
=&\,
\be(x)_{13}. 
\end{align*}
Hence $\be(x)\in \RG\oti L^\infty(\bH\setm\bG)$. 
By weak continuity of $\be$, 
$\be(L^\infty(\bH\setm\bG))\subs\RG\oti L^\infty(\bH\setm\bG)$. 
\end{proof}

In fact assuming the co-amenability $\bG$, 
We prove the converse statement of the previous lemma 
in Theorem \ref{thm: exco}. 

Let $B\subs\lG$ be a right coideal. 
We denote by $L^2(B)$ 
the norm closure of the space $B\hat{1}_h$. 
Assume that $B$ has the expectation property, that is, 
there exists a conditional expectation 
$E_B\col \lG\ra B$ such that $E_B$ preserves the Haar state $h$. 
Define the Jones projection $e_B\col \ltG\ra L^2(B)$ 
by 
\[
e_B (x \hat{1}_h)= E_B(x)\hat{1}_h 
\quad\mbox{for}\ x\in \lG.
\] 
Let $B\subs \lG\subs \lG\vee \{e_B \}''=:B_1$ 
be the basic extension. 
The main properties of $e_B$ are as follows 
(see \cite[Lemma 3.2]{Ko1} and \cite[p. 312]{Ta}). 

\begin{lem}
With the above settings, one has 
\begin{enumerate}

\item 
$e_B x e_B=E_B(x)e_B$ for $x\in\lG$, 

\item 
$B=\lG\cap \{e_B\}'$, 

\item 
$Je_B=e_B J$, $\vDe_h^{it}e_B=e_B\vDe_h^{it}$ for all $t\in\R$, 

\item 
$B_1=J B'J$. 

\end{enumerate}
\end{lem} 

Set $\hB=B'\cap \LG$. 
Then $\hB$ is a left coideal of $\LG$ as is shown 
in Lemma \ref{lem: lattice}. 
The following lemma is proved in 
\cite[Theorem 4.6]{ILP}, which treats the Kac algebra case. 
The proof can be adapted to the quantum group case. 

\begin{lem}\label{lem: Jones}
With the above setting, one has 
\begin{enumerate}

\item $e_B\in \hB$, 

\item $\de(E_B(x))=(E_B\oti\id)(\de(x))$ for all $x\in\lG$, 

\item $\hB'\cap \lG=B$. 

\end{enumerate}
\end{lem}
\begin{proof}
(1). 
Since $B$ is a right coideal, 
we see that 
$V_\bG(L^2(B)\oti \ltG)\subs L^2(B)\oti \ltG$. 
Hence $V_\bG(e_B\oti1)=(e_B\oti1)V_\bG(e_B\oti1)$. 
Since $(J\oti \hJ)V_\bG(J\oti \hJ)=V_\bG^*$ 
and $Je_B J=e_B$, 
$V_\bG(e_B\oti1)=(e_B\oti1)V_\bG$. 
Hence $e_B\in B'\cap \LG=\hB$. 

(2). 
Let $x\in \lG$. Since $e_B\in\LG$ and $V_\bG\in \RG\oti\lG$, 
we have
\begin{align*}
(E_B\oti\id)(\de(x))(e_B\oti1)
=&\,
(e_B\oti1)\de(x)(e_B\oti1)
=
(e_B\oti1)V_\bG(x\oti1)V_\bG^*(e_B\oti1)
\\
=&\,
V_\bG (e_B\oti1)(x\oti1)(e_B\oti1)V_\bG^*
=
V_\bG (E_B(x)e_B\oti1)V_\bG^*
\\
=&\,
V_\bG (E_B(x)\oti1)V_\bG^*(e_B\oti1)
\\
=&\,
\de(E_B(x))(e_B\oti1). 
\end{align*}
In particular, we have 
$\de(E_B(x))(\hat{1}_h\oti\hat{1}_h)
=(E_B\oti\id)(\de(x))(\hat{1}_h\oti\hat{1}_h)$. 
Since $\hat{1}_h\oti \hat{1}_h$ is a separating vector 
for $\lG\oti\lG$, 
we have $\de(E_B(x))=(E_B\oti\id)(\de(x))$. 

(3). It follows from Lemma \ref{lem: lattice} or the direct argument 
as follows. 
It is clear that $B\subs \hB'\cap \lG$. 
Since $e_B\in\hB$ and $B=\{e_B\}'\cap \lG$, 
$B=\hB'\cap \lG$. 
\end{proof}

Next we consider a relation between $\hB$ and $B_1$. 
Let $\al$ be the right $\bG$-action $\al$ on $B(\ltG)$ 
defined in \S \ref{subsect: al-be}. 
Since $e_B\oti1$ commutes with $V_\bG$, 
$B_1$ is globally invariant under the right action $\al$. 
Let $X$ be a globally invariant subspace in $B_1$. 
The set of the fixed point elements of $X$ under $\al$ 
is denoted by $X^\bG$. 
Note that $X^\bG=X\cap\LG$ holds. 

\begin{lem}\label{lem: eB}
The following equalities hold:
\begin{enumerate}

\item $\hR(e_B)=e_B$. 

\item $(e_B\oti1)W_\bG^*(1\oti e_B)=W_\bG^*(e_B\oti e_B)$. 

\item $\wdh{B}=J B_1^\bG J=\ovl{\C+J(\lG e_B\lG)^\bG J}^{\,\rm{w}}$. 

\end{enumerate}
\end{lem}
\begin{proof}
(1). 
Since $e_B^*=e_B$ and $Je_B=e_B J$, we have $\hR(e_B)=Je_B^* J=e_B$. 

(2). 
Take any $x,y\in\lG$. Then we have 
\begin{align*}
(e_B\oti1)W_\bG^*(1\oti e_B) (x\hat{1}_h\oti y\hat{1}_h)
=&\,
(e_B\oti1)W_\bG^* (x\hat{1}_h\oti E_B(y)\hat{1}_h)
\\
=&\,
(e_B\oti1) (\de(E_B(y)) x\hat{1}_h\oti \hat{1}_h)
\\
=&\,
\de(E_B(y)) (E_B(x)\hat{1}_h\oti \hat{1}_h)
\\
=&\,
W_\bG^*(e_B\oti e_B)(x\hat{1}_h\oti y\hat{1}_h). 
\end{align*}

Hence the desired equality holds. 

(3). 
Since $B_1=JB'J$ and $J\LG J=\LG$, we have 
\begin{align*}
J\wdh{B}J
=&\,
J(B'\cap\LG)J
=
JB'J\cap J\LG J
\\
=&\,
B_1\cap \LG
=
B_1^\bG. 
\end{align*}

By Lemma \ref{lem: Jones}, the $*$-subalgebra $\lG+\lG e_B \lG$ 
is weakly dense in $B_1$. 
Hence we have $B_1^\bG=\ovl{\C + (\lG e_B\lG)^\bG}^{\,\rm{w}}$. 
\end{proof} 

\begin{lem}\label{lem: AGeB}
One has 
\[
(\AG e_B \AG)^\bG
=\{(\id\oti\om)(\De(e_B))\mid \om\in \LG_*^{\rm{fin}}\}. 
\]
\end{lem}
\begin{proof}
Let $\pi,\si\in\IG$, $i,j\in I_\pi$ and $k,\el\in I_\rho$. 
Then 
\begin{align*}
E_\al(v_{\pi_{i,j}}e_B v_{\rho_{k,\el}}^*)
=&\,
(\id\oti h)(\al(v_{\pi_{i,j}}e_B v_{\rho_{k,\el}}^*))
\\
=&\,
\sum_{m\in I_\pi}\sum_{n\in I_\rho}
v_{\pi_{i,m}}e_B v_{\rho_{k,n}}^*
h(v_{\pi_{m,j}}v_{\rho_{n,\el}}^*)
\\
=&\,
\sum_{m\in I_\pi}\sum_{n\in I_\rho}
\de_{\pi,\rho}\de_{m,n}\de_{j,\el}D_\pi^{-1} F_{\pi_{j,j}}
v_{\pi_{i,m}}e_B v_{\rho_{k,n}}^*
\\
=&\,
D_\pi^{-1} F_{\pi_{j,j}}
\de_{\pi,\rho}\de_{j,\el}
\sum_{m\in I_\pi}
v_{\pi_{i,m}}e_B v_{\pi_{k,m}}^*.
\end{align*}
Hence 
\[(\AG e_B\AG)^\bG 
=\spa\Big{\{}\sum_{m\in I_\pi}v_{\pi_{i,m}}e_B v_{\pi_{k,m}}^*
\,\Big{|}\,
i,k\in I_\pi, \pi\in\IG\Big{\}}.
\] 
Using the description of $W_\bG$ in (\ref{eq: decomp}), 
we have 
\[
\De(e_B)(1\oti 1_\pi)
=
\sum_{i,j,k\in I_\pi}
v_{\pi_{i,j}}e_B v_{\pi_{k,j}}^*\oti f_{\pi_{i,k}}. 
\]
Hence 
\[
(\AG e_B \AG)^\bG
=\{(\id\oti\om)(\De(e_B))\mid \om\in \LG_*^{\rm{fin}}\}. 
\]
\end{proof}

\begin{lem}
One has $B_1^\bG=\ovl{(\lG e_B \lG)^\bG}^{\,\rm{w}}$. 
\end{lem}
\begin{proof}
By Lemma \ref{lem: eB}, 
it suffices to show that the unit of $B_1$ is equal to 
that of $\ovl{(\lG e_B \lG)^\bG}^{\,\rm{w}}$. 
Let $p\in B_1^\bG$ be the unit of the von Neumann subalgebra 
generated by $(\lG e_B \lG)^\bG$. 
Set $q=1-p$. 
Then by the previous lemma, we have $\De(e_B)(q\oti1)=0$. 
Since $\hR(e_B)=e_B$ by Lemma \ref{lem: eB}, 
we have 
\begin{align*}
(1\oti \hR(q))\De(e_B)
=&\,
(1\oti \hR(q))\De(\hR(e_B))
=
(1\oti \hR(q)) (\hR\oti\hR)(\De^\op(e_B))
\\
=&\,
(\hR\oti\hR)(\De^\op(e_B)(1\oti q))
=0. 
\end{align*}
The equality 
$\De(e_B)=W_\bG(e_B\oti1)W_\bG^*$ 
yields $(1\oti \hR(q))W_\bG(e_B\oti1)=0$, 
and 
$(E_B\oti\id)(W_\bG^*(1\oti \hR(q))W_\bG)=0$. 
Since $E_B$ is faithful, 
we have $q=0$. 
\end{proof}

\begin{lem}\label{lem: B_1} One has
\[
B_1^\bG=\ovl{\{(\id\oti\om)(\De(e_B))\mid \om\in \LG_*\}}^{\,\rm{w}},
\quad
\hB=\ovl{\{(\om\oti\id)(\De(e_B))\mid \om\in \LG_*\}}^{\,\rm{w}}.
\]
\end{lem}
\begin{proof}
The left equality follows from Lemma \ref{lem: AGeB} 
and the previous lemma. 
Take $\om\in\LG_*$. 
By (\ref{eq: DehR}) and Lemma \ref{lem: eB}, we have
\[
J(\id\oti\om)(\De(e_B))J
=(\id\oti\ovl{\om}\circ \hR)((\hR\oti\hR)(\De(e_B)))
=(\ovl{\om}\circ \hR\oti\id)(\De(e_B))
\]
Hence the right equality holds since $\hB=J B_1^\bG J$. 
\end{proof}

To construct a left invariant weight on $\hB$, 
we make use of theory of spatial derivatives and 
operator valued weights which have been introduced 
in \cite{Co}, \cite{Ha1} and \cite{Ha2}. 
Readable explanations on them are presented in \cite{Ko2} 
and we freely use the notations there. 
Let $E_B^{-1}\col B'\ra \lG'$ be the operator valued weights 
associated with $E_B\col \lG\ra B$. 
It is characterized by the following equality 
on spatial derivatives:
\begin{equation}\label{eq: spatial}
\frac{d\om'\circ E_B^{-1}}{d\om}
=
\frac{d\om'}{d\om\circ E_B},
\end{equation}
where $\om'$ and $\om$ are faithful normal semifinite 
weights on $\lG'$ and $B$, respectively. 
The equality $E_B^{-1}(e_B)=1$ holds \cite[Lemma 3.1]{Ko1}. 
We define the $*$-subalgebra $\hB_{0}$ of $\hB$ by 
\begin{equation}\label{eq: hB0}
\hB_0=J(\AG e_B\AG)^\bG J
=\{(\om\oti\id)(\De(e_B))\mid \om\in \LG_*^{\rm{fin}}\}.
\end{equation}

\begin{lem}\label{lem: semi}
The operator valued weight $E_B^{-1}$ is finite on $\hB_0$. 
In particular, $E_B^{-1}$ is semifinite on $\hB$. 
\end{lem}
\begin{proof}
Set a unitary $W_\bhG=(J\oti J)(W_\bG)_{21}^*(J\oti J)\in \LG\oti \lG'$. 
Then for any $x\in \LG$, we have 
\begin{align*}
W_\bhG^*(1\oti x)W_\bhG
=&\,
(J\oti J) (W_\bG)_{21}(1\oti JxJ)(W_\bG)_{21}^* (J\oti J)
\\
=&\,
(J\oti J) (W_\bG)_{21}(1\oti \hR(x^*))(W_\bG)_{21}^* (J\oti J)
\\
=&\,
(J\oti J) \De(\hR(x^*))_{21} (J\oti J)
\\
=&\,
(\hR\oti\hR)(\De(\hR(x))_{21})
\\
=&\,
\De(x). 
\end{align*}

Now take a positive functional $\om\in \LG_*^{\rm{fin}}$. 
Then by $\lG'$-linearity of $E_B^{-1}$, 
we have 
\begin{align*}
E_B^{-1}((\om\oti\id)(\De(e_B)))
=&\,
E_B^{-1}((\om\oti\id)(W_\bhG^*(1\oti e_B)W_\bhG))
\\
=&\,
(\om\oti\id)(W_\bhG^*(1\oti E_B^{-1}(e_B))W_\bhG)
\\
=&\,
\om(1)<\infty.
\end{align*}
Hence $E_B^{-1}$ is finite on $\hB_0$. 
\end{proof}

Define a state $h'\in \lG'$ 
by $h'(x)=h(Jx^*J)$ for $x\in \lG'$. 
Set a faithful normal semifinite weight 
$\vph=h'\circ E_B^{-1}$ on $B'$. 
We denote by $\vph_\hB$ the restriction $\vph$ on $\hB$. 

\begin{lem}\label{lem: preserve}
The modular automorphism group $\{\si_t^{\vph_\hB}\}_{t\in\R}$ preserves 
$\hB_{0}$. 
\end{lem}
\begin{proof}
By (\ref{eq: spatial}), 
we obtain
\[
\frac{d\vph}{d h_B}
=
\frac{d h'\circ E_B^{-1}}{d h|_B}
=
\frac{d h'}{ d h|_B\circ E_B}
=
\frac{d h'}{ d h}
=
\vDe_h^{-1}. 
\]
Hence $\si_t^{\vph}=\Ad \vDe_h^{-it}$ on $B'$. 
Since $\Ad \vDe_h^{-it}=\hta_{t}$ on $\LG$ 
by (\ref{eq: vDP}) and (\ref{eq: hRhta}), 
the von Neumann subalgebra $\hB\subs B'$ 
is globally invariant under $\si_t^{\vph}$. 
Hence 
$\si_t^{\vph_\hB}=\si_t^{\vph}|_\hB=\hta_{t}|_\hB$ 
by uniqueness of the modular automorphism group. 
Since $e_B$ commutes with $\vDe_h^{-it}$, 
$\hta_{t}(e_B)=\vDe_h^{-it} e_B\vDe_h^{it}=e_B$. 
Then for any $\om\in\LG_*^{\rm{fin}}$, 
\[
\si_t^{\vph_\hB}((\om\oti\id)(\De(e_B)))
=
\hta_{t}((\om\oti\id)(\De(e_B)))
=
(\om\circ\hta_{-t}\oti\id)(\De(e_B))
\]
by (\ref{eq: DehR}) and $\hta_t(e_B)=e_B$. 
By definition of $\hB_0$ (\ref{eq: hB0}), 
we have $\si_t^{\vph_\hB}(\hB_0)\subs\hB_0$. 
\end{proof}

By previous lemma, 
the following push-down lemma is proved 
in a similar way to \cite[Proposition 2.2]{ILP}. 

\begin{lem}\label{lem: push-down}
For any $x\in n_{\vph_\hB}$, $e_B x=e_B E_B^{-1}(e_B x)$ holds. 
\end{lem}

\begin{lem}\label{lem: B0}
One has 
$(\AG e_B \lG)^\bG=(\AG e_B \AG)^\bG$. 
\end{lem}
\begin{proof}
It suffices to show 
$(\lG_\pi e_B \lG)^\bG\subs \lG_\pi e_B \lG_\pi^*$, 
where $\lG_\pi$ is the linear span of $\{v_{\pi_{i,j}}\}_{i,j\in I_\pi}$. 
Take any $x\in \lG$. Then 
\[
E_\al(v_{\pi_{i,j}} e_B x)
=\sum_{k\in I_\pi} 
v_{\pi_{i,k}}e_B(\id\oti h)((1\oti v_{\pi_{k,j}})\de(x))
\in \lG_\pi e_B \lG_\pi^*, 
\]
because $(\id\oti h)((1\oti v_{\pi_{k,j}})\de(x))\in \lG_\pi^*$ 
for any $x\in\lG$.  
\end{proof}

\begin{lem}
The weight $\vph_\hB$ is left invariant, that is, 
for any $x\in m_{\vph_\hB}$ and $\om\in\LG_*$, 
one has 
\[
(\om\oti\vph_\hB)(\De(x))=\om(1)\vph_\hB(x).
\]
\end{lem}
\begin{proof}
Take $\th\in\LG_*^{\rm{fin}}$ and 
set $x=(\th\oti\id)(\De(e_B))\in \hB_0$. 
By using the proof of Lemma \ref{lem: semi}, 
we have $E_B^{-1}(x)=\th(1)$. 
Then for any $\om\in\LG_*$, 
\begin{align*}
(\om\oti\vph_\hB)(\De(x))
=&\,
(\om\oti h'\circ E_B^{-1})(W_\bhG^*(1\oti x)W_\bhG)
\\
=&\,
(\om\oti h')(W_\bhG^*(1\oti E_B^{-1}(x))W_\bhG)
\\
=&\,
\om(1)\th(1)
\\
=&\,
\om(1)\vph_\hB(x). 
\end{align*}
Hence on $\hB_{0}$, the left invariance holds. 

Take a bounded sequence $\{u_n\}_{n\in \N}$ in $\hB_0$ 
which strongly converges to $1$. 
Let $y\in n_{\vph_\hB}$. 
Then by Lemma \ref{lem: push-down}, 
$u_n y \in J (\AG e_B \lG)J$ for each $n\in \N$. 
Since $J u_n y J\in J\hB J=B_1^\bG$, we see that 
$u_n y \in J (\AG e_B \lG)^\bG J$. 
Using Lemma \ref{lem: B0}, 
we see that $u_n y \in J (\AG e_B \AG)^\bG J=\hB_0$. 
Hence for any positive $\om\in \LG_*$, 
we have 
\begin{equation}\label{eq: yuuy}
(\om\oti\vph_\hB)(\De(y^* u_n^* u_n y))=\om(1)\vph_\hB(y^* u_n^* u_n y). 
\end{equation}
Taking the limit as $n\to\infty$ in (\ref{eq: yuuy}), 
we have 
\begin{align*}
(\om\oti\vph_\hB)(\De(y^* y))
\leq&\,
\lim_{n\to\infty}
(\om\oti\vph_\hB)(\De(y^* u_n^* u_n y))
\\
=&\,
\lim_{n\to\infty}
\om(1)\vph_\hB(y^* u_n^* u_n y)
\\
=&\,
\om(1)\vph_\hB(y^* y)<\infty. 
\end{align*}
In particular, 
the map $x\in\hB\mapsto (\om\oti\vph_\hB)(\De(y^* x y))$ 
is a normal functional. 
Hence again taking the limit as $n\to\infty$ in (\ref{eq: yuuy}), 
we have 
\[
(\om\oti\vph_\hB)(\De(y^* y))=\om(1)\vph_\hB(y^* y). 
\]
Therefore $\vph_\hB$ is left invariant on $m_{\vph_\hB}$. 
\end{proof}

We summarize our arguments as follows. 
\begin{thm}
Let $B\subs \lG$ be a right coideal. 
If $B$ has the expectation property, 
then the left coideal $\hB$ has a left invariant 
faithful normal semifinite weight for the left action 
of $\LG$. 
\end{thm}

Next we study a right coideal of $\bG$ endowed with 
the coaction symmetry. 

\begin{lem}
Let $B\subs \lG$ be a right coideal. 
Then $B$ has the coaction symmetry 
if and only if 
$\De(\wdh{B})\subs\wdh{B}\oti\wdh{B}$. 
\end{lem}
\begin{proof}
We know $\De(\hB)\subs \LG\oti \hB$. 
Let $x\in B$ and $y\in\wdh{B}$. 
Then we have 
\begin{align*}
\De(y)(x\oti1)
=&\,
W_\bG (y\oti1) W_\bG^* (x\oti1)
\\
=&\,
(1\oti U_\bG)(V_\bG)_{21}(1\oti U_\bG)
(y\oti1)
(1\oti U_\bG)(V_\bG)_{21}^*(1\oti U_\bG)(x\oti1)
\\
=&\,
(1\oti U_\bG)(V_\bG)_{21}
(y\oti1)
\be(x)_{21}(V_\bG)_{21}^*(1\oti U_\bG)
\end{align*}
and similarly 
\[
(x\oti1)\De(y)
=
(1\oti U_\bG)(V_\bG)_{21}
\be(x)_{21}(y\oti1)
(V_\bG)_{21}^*(1\oti U_\bG).
\]
Hence $\De(\wdh{B})\subs \wdh{B}\oti\wdh{B}$ if and only if 
$\be(B)\subs \RG\oti (\wdh{B}'\cap \lG)=\RG\oti B$. 
\end{proof}

Set a weight $\ps_\hB=\vph_\hB\circ\hR$ on $\hB$. 
Using (\ref{eq: DehR}), 
we see that $\ps_\hB$ is right invariant. 
Therefore $(\hB,\De)$ is a quantum group 
in the sense of \cite{KV}. 
Clearly the restriction $\hvep|_\hB$ is a normal counit on $\hB$. 
Hence $(\hB,\De)$ is a discrete quantum group. 
The counit $\hvep|_\hB$ is given 
by cutting elements at the Jones projection $e_B$ as follows. 

\begin{lem}
The equality $x e_B=\hvep(x)e_B=e_B x$ holds for all $x\in \hB$. 
In particular, 
the Jones projection 
$e_B$ is a minimal central projection of $\hB$.
\end{lem}
\begin{proof}
By Lemma \ref{lem: eB}, 
$(e_B\oti1)W_\bG^*(1\oti e_B)=W_\bG^* (e_B\oti e_B)$. 
Then we have 
\[
\De(e_B)(1\oti e_B)
=
W_\bG(e_B\oti1)W_\bG^*(1\oti e_B)
=e_B\oti e_B.
\]
Taking the adjoint of the above equality, we have 
\[
\De(e_B)(1\oti e_B)=e_B\oti e_B=(1\oti e_B)\De(e_B). 
\] 
Take any $\om\in\LG_*$ and set $x=(\om\oti\id)(\De(e_B))$. 
Then we have 
\begin{align*}
x e_B
=&\,
(\om\oti\id)(\De(e_B))e_B
=
(\om\oti\id)(\De(e_B)(1\oti e_B))
\\
=&\,
(\om\oti\id)(e_B\oti e_B)
=
\om(e_B) e_B
=
\hvep(x)e_B. 
\end{align*} 
Similarly we obtain $e_B x=\hvep(x)e_B$. 
Then the desired equality holds for all $x\in \hB$ 
by Lemma \ref{lem: B_1}. 
\end{proof}

We summarize our arguments as follows. 

\begin{thm}\label{thm: discrete}
Let $B\subs \lG$ be a right coideal. 
If $B$ has the expectation property and the coaction symmetry, 
then the pair $(\hB,\De)$ is a discrete quantum group. 
\end{thm}

Under assumption on co-amenability, 
we obtain the following characterization of 
right coideals of quotient type. 

\begin{thm}\label{thm: exco}
Let $\bG$ be a co-amenable compact quantum group 
and $B\subs \lG$ a right coideal. 
Then 
$B$ is of quotient type if and only if 
$B$ has the expectation property and the coaction symmetry. 
\end{thm}
\begin{proof}
We have already proved the ``only if'' part in Lemma \ref{lem: exco}. 
So, it suffices to show the ``if'' part. 
By the previous theorem, 
$(\hB,\De)$ is a discrete quantum group. 
Let $\bH=(\CH,\de_\bH)$ be a compact quantum group 
such that the discrete quantum groups 
$(\LH,\De_\bH)$ and $(\hB,\De)$ 
are isomorphic. 
We identify $\hB$ with $\LH$. 
By Lemma \ref{lem: algsub-corresp}, 
$\bH$ is represented as an algebraic quantum subgroup of $\bG$, 
that is, 
there exists a restriction map 
$r_\bH\col A(\bG)\ra A(\bH)$ 
such that $r_\bH^t\col \hB\ra\LG$ is a given inclusion. 
Moreover by Lemma \ref{lem: algsub}, 
$\{\bH,r_\bH\}$ is in fact a quantum subgroup of $\bG$. 
By Lemma \ref{lem: above} and Lemma \ref{lem: Jones}, we have 
$B=\hB'\cap \lG=\LH'\cap \lG=L^\infty(\bH\setm\bG)$. 
\end{proof}

\section{Application to classification of Poisson boundaries}

In the rest of this paper, we determine the Poisson boundary 
for a co-amenable compact quantum group with the commutative 
fusion rules. 
We also compute the Poisson boundary for a 
$q$-deformed classical compact Lie group. 

\subsection{Poisson boundaries.}

We briefly recall the notion of the Poisson boundary for 
a discrete quantum group. 
We refer to \cite{I1}, \cite{INT} and \cite{NT} 
for definitions of terminology. 

Let $\ph_\pi\in B(H_\pi)_*$ be the right $\bG$-invariant state. 
Define a transition operator $P_\pi$ on $\RG$ by 
$P_\pi(x)=(\id\oti \ph_\pi)(\De_R(x))$ 
for $x\in \RG$. 
For a probability measure $\mu$ on $\IG$, we set 
a non-commutative Markov operator 
\[
P_\mu=\sum_{\pi\in\IG}\mu(\pi)P_\pi. 
\]
Then for a generating measure $\mu$, we define an 
operator system 
\[
H^\infty(\bhG,P_\mu)=\{x\in \RG\mid P_\mu(x)=x\}. 
\]
It has the von Neumann algebra structure defined by 
\[
x\cdot y= \lim_{n\to\infty} P_\mu^n(xy) 
\quad \mbox{for}\  x,y\in H^\infty(\bhG,P_\mu), 
\]
where the limit is taken in the strong topology \cite[Theorem 3.6]{I1}. 
The von Neumann algebra $H^\infty(\bhG,P_\mu)$ is called 
the (non-commutative) Poisson boundary of $\{\RG,P_\mu\}$. 
Assuming that the fusion algebra of $\bG$ is commutative, 
we know 
that the Poisson boundary does not depend on the generating measure 
\cite[Proposition 1.1]{INT}, that is, we have 
\[
H^\infty(\bhG,P_\mu)=\{x\in\RG\mid P_\pi(x)=x\ \mbox{for all}\ \pi\in\IG\}. 
\]
We are interested only in such a compact quantum group, 
so we write simply $H^\infty(\bhG)$ for $H^\infty(\bhG,P_\mu)$. 
The Poisson integral $\Th\col \lG\ra\RG$ is defined by 
\begin{equation}\label{eq: Poisson}
\Th=(\id\oti h)\circ\be, 
\end{equation}
which maps $\lG$ into $H^\infty(\bhG)$ \cite[Lemma 3.8]{I1}. 
It is $\bhG$-$\bG$-equivariant 
in the following sense:
\[
\al\circ \Th=(\Th\oti\id)\circ\de,
\quad
\De_R\circ\Th=(\id\oti\Th)\circ\be.
\]

Now we consider how the inverse map of $\Th$ can be constructed. 
That will be written as a similar form to (\ref{eq: Poisson}), 
that is, $\al$ and some state $\om$ on $\RG$ 
will take the place of $\be$ and $h$, respectively. 
Then we consider the map $\RG\ni x\mapsto (\om\oti\id)(\al(x))\in\lG$. 
The inverse of $\Th$ (if it exists) 
should be $\bhG$-$\bG$-equivariant as $\Th$ is.
Here, we have to realize what property of $h$ 
derives the bi-equivariance of $\Th$ in the proof of \cite[Lemma 3.8]{I1}. 
While the $\bhG$-equivariance follows by definition of the left action, 
the $\bG$-equivariance does by right invariance of $h$. 
Hence the state $\om$ has to satisfy the left invariance 
for the coproduct $\De_R$. 

From now we assume the amenability of $(\RG,\De_R)$. 
Let $m\in \RG^*$ be a left invariant mean. 
Although $m$ is non-normal in general, 
we can consider a unital completely positive map 
$m\oti\id\col \RG\oti M\ra M$ 
for any von Neumann algebra $M$. 
Indeed for any $x\in \RG\oti M$, we define 
an element $(m\oti\id)(x)\in M=(M_*)^*$ by 
$\om((m\oti\id)(x))=m((\id\oti\om)(x))$ 
for $\om\in M_*$. 
Let $M$ and $N$ be von Neumann algebras 
and $T\col M\ra N$ a normal completely bounded map. 
Then 
$(m\oti\id_N)\circ(\id\oti T)=T\circ (m\oti\id_M)$ holds. 
Here we need to assume the normalcy of $T$. 
In particular, we have 
$(m\oti\id)((1\oti a)x (1\oti b))=a (m\oti\id)(x)b$ 
for all $x\in \RG\oti M$ and $a,b\in M$. 

Now define a unital completely positive map 
$\La\col \RG\ra \lG$ by 
\[
\La(x)=(m\oti\id)(\al(x))
\quad\mbox{for}\ x\in \RG. 
\]
It turns out that $\La$ is normal 
on the Poisson boundary $H^\infty(\bhG)$ in the next lemma. 
Moreover by Theorem \ref{thm: poisson}, 
we will see that $\La$ is actually equal to $\Th^*$ 
defined in \cite{INT}. 
This means that $\La$ does not depend 
on the choice of an left invariant mean. 

\begin{lem}\label{lem: equiv}
The following statements hold. 
\begin{enumerate}

\item 
$h\circ\La=\hvep$ on $H^\infty(\bhG)$. 
In particular, 
$\La$ is a faithful normal map. 

\item 
$\La$ is $\bhG$-$\bG$-equivariant, that is, 
\[
\de\circ\La=(\La\oti\id)\circ\al, \quad
\be\circ\La=(\id\oti\La)\circ\De_R. 
\]
\end{enumerate}
\end{lem}
\begin{proof}
(1). 
Since the fusion algebra of $\bG$ is commutative, 
the action $\al$ on $H^\infty(\bhG)$ is ergodic \cite[Corollary 3.7]{I1}. 
Hence $E_\al(x)=\hvep(x)1$ for all $x\in H^\infty(\bhG)$. 
Then $h(\La(x))=m(E_\al(x))=m(\hvep(x)1)=\hvep(x)$. 
Since $h$ and $\hvep$ are faithful normal states \cite[Theorem 3.6]{I1}, 
we conclude that $\La$ is a faithful normal map. 

(2). 
First we show the $\bG$-equivariance of $\La$. 
Let $x\in \RG$. Then we have 
\begin{align*}
\de(\La(x))
=&\,
\de((m\oti\id)(\al(x)))
=
(m\oti\id\oti\id)((\id\oti\de)(\al(x)))
\\
=&\,
(m\oti\id\oti\id)((\al\oti\id)(\al(x)))
\\
=&\,
((m\oti\id)\circ\al\oti\id)(\al(x))
=
(\La\oti\id)(\al(x)). 
\end{align*}

Next we show the $\bhG$-equivariance of $\La$. 
The left invariance of $m$ yields $(\id\oti m)(\De_R(y))=m(y)1$ 
for all $y\in\RG$. 
Take any $x\in \RG$ and then 
\begin{align*}
&(\id\oti\La)(\De_R(x))
\\
=&\,
(\id\oti m\oti\id)
((\id\oti \al)(\De_R(x)))
\\
=&\,
(\id\oti m\oti\id)
((V_\bG)_{23}(V_\bG)_{12}^* 
(1\oti x\oti 1)(V_\bG)_{12}(V_\bG)_{23}^*)
\\
=&\,
(\id\oti m\oti\id)
((V_\bG)_{13}^*(V_\bG)_{12}^*(V_\bG)_{23}
(1\oti x\oti 1)(V_\bG)_{23}^*(V_\bG)_{12}(V_\bG)_{13})
\\
=&\,
V_\bG^* 
(\id\oti m\oti\id)((\De_R\oti\id)(\al(x))) 
V_\bG
\\
=&\,
V_\bG^* 
((\id\oti m)\circ\De_R\oti\id)(\al(x)) 
V_\bG
\\
=&\,
V_\bG^* 
(1\oti (m\oti\id)(\al(x)))
V_\bG
=
\be(\La(x)). 
\end{align*}
\end{proof}

\begin{lem}\label{lem: ThLa}
For any $x\in H^\infty(\bhG)$, 
$\Th(\La(x))=x$. 
In particular, $\La$ is a faithful normal 
$*$-homomorphism. 
\end{lem}
\begin{proof}
Let $x\in H^\infty(\bhG)$. 
By Lemma \ref{lem: equiv}, we have 
\begin{align*}
\Th(\La(x))
=&\,
(\id\oti h)(\be(\La(x)))
=
(\id\oti h)((\id\oti\La)(\De_R(x)))
\\
=&\,
(\id\oti \hvep)(\De_R(x))
=
x. 
\end{align*}
The multiplicativity of $\La$ is shown as follows. 
For any $x\in H^\infty(\bhG)$, we have 
\begin{align*}
x^* \cdot x
=&\,
\Th(\La(x))^*\cdot \Th(\La(x))
\leq
\Th(\La(x)^* \La(x))
\\
\leq&\,
\Th(\La(x^* \cdot x))
=
x^* \cdot x. 
\end{align*}
Hence we obtain $\Th(\La(x)^* \La(x))=\Th(\La(x^*\cdot x))$, 
namely, $\La(x)^* \La(x)=\La(x^*\cdot x)$. 
This implies that $\La$ is a $*$-homomorphism. 
\end{proof}

We determine the multiplicative domain of $\Th$. 

\begin{lem}\label{lem: prepoisson}
Set $B=\La(H^\infty(\bhG))$. 

\begin{enumerate}

\item
$B$ is a right coideal 
with the expectation property and the coaction symmetry. 

\item 
The conditional expectation onto $B$ is $\bhG$-equivariant. 

\item
$B$ coincides with 
the multiplicative domain of $\Th$. 
Moreover, the Poisson integral 
$\Th\col B\ra H^\infty(\bhG)$ is an isomorphism.

\end{enumerate}
\end{lem}
\begin{proof}
(1). 
By the previous lemma, $B\subs\lG$ is a von Neumann subalgebra. 
Since $\La$ is $\bhG$-$\bG$-equivariant by Lemma \ref{lem: equiv}, 
$B$ is globally invariant under the actions $\al$ and $\be$. 
Hence $B$ is a right coideal with the coaction symmetry. 
Set $E_B=\La\circ\Th$. 
Then by the previous lemma, 
we see that 
$E_B$ is a faithful normal conditional expectation onto $B$. 
For $x\in\lG$, we have 
\[
h(E_B(x))
=
h(\La(\Th(x)))
=
\hvep(\Th(x))
=
h(x).  
\]
Hence $B$ has the expectation property. 

(2). 
Since $\Th$ and $\La$ are $\bhG$-equivariant, 
$E_B=\Th\circ\La$ is also $\bhG$-equivariant. 

(3). 
Let $\mD\subs \lG$ be the multiplicative domain of $\Th$. 
It is easy to see that $B\subs \mD$. 
We show the converse inclusion. 
Let $x\in\mD$. 
Then by definition, 
we have $\Th(x^* x)=\Th(x)^*\cdot \Th(x)$ 
and $\Th(x x^*)=\Th(x)\cdot \Th(x)^*$. 
Applying $\La$ to both the sides of the equalities, 
we have 
$E_B(x^* x)=E_B(x)^* E_B(x)$ 
and $E_B(xx^*)=E_B(x)E_B(x)^*$. 
This immediately yields that $x\in B$, and $B=\mD$. 
Hence $\Th\col B\ra H^\infty(\bhG)$ 
is a faithful normal $*$-homomorphism. 
The surjectivity of $\Th$ follows from 
$\Th(B)=\Th(\La(H^\infty(\bhG)))=H^\infty(\bhG)$. 
\end{proof}

\begin{lem}\label{lem: Kac}
Let $\bH$ be a quantum subgroup of $\bG$. 
Then $\bH$ is of Kac type 
if and only if 
the expectation $E_\bH\col \lG\ra L^\infty(\bH\setm\bG)$ 
is $\bhG$-equivariant, that is, 
it satisfies 
$(\id\oti E_\bH)\circ\be=\be\circ E_\bH$. 
\end{lem}
\begin{proof}
For any $x\in \CG$, we have 
\begin{align*}
(\id\oti E_\bH)(\be(x))
=&\,
(\id\oti h_\bH\circ r_\bH \oti\id)
((\id\oti\de)(V_\bG^* (1\oti x)V_\bG))
\\
=&\,
(\id\oti h_\bH\circ r_\bH\oti\id)
((V_\bG)_{13}^*(V_\bG)_{12}^*(1\oti \de(x))(V_\bG)_{12}(V_\bG)_{13})
\\
=&\,
V_\bG^* (\id\oti h_\bH\oti\id)
((V_\bH)_{12}^*(1\oti(r_\bH\oti\id)(\de(x)))(V_\bH)_{12}) 
V_\bG.
\end{align*}
This is equal to $\be(E_\bH(x))=V_\bG^*(1\oti E_\bH(x))V_\bG$ 
if and only if  
\[
(\id\oti h_\bH\oti\id)
((V_\bH)_{12}^*(1\oti (r_\bH\oti\id)(\de(x)))(V_\bH)_{12}) 
=
1\oti (h_\bH\oti\id)((r_\bH\oti\id)(\de(x))). 
\]
Multiplying $1\oti1\oti y$, $y\in C(\bG)$ from the right, 
we have 
\begin{align*}
&(\id\oti h_\bH\oti\id)
((V_\bH)_{12}^*(1\oti (r_\bH\oti\id)(\de(x)(1\oti y)))(V_\bH)_{12}) 
\\
=&\,
1\oti (h_\bH\oti\id)((r_\bH\oti\id)(\de(x)(1\oti y))). 
\end{align*}
Since the subspace $\de(C(\bG))(\C\oti C(\bG))$ 
is dense in $C(\bG)\oti C(\bG)$, 
we have for all $x\in C(\bG)\oti C(\bG)$, 
\[
(\id\oti h_\bH\oti\id)
((V_\bH)_{12}^*(1\oti (r_\bH\oti\id)(x))(V_\bH)_{12}) 
=
1\oti (h_\bH\oti\id)((r_\bH\oti\id)(x))). 
\]
More precisely, we have 
\[
(\id\oti h_\bH)
(V_\bH^*(1\oti z) V_\bH) 
=
h_\bH(z)1 
\quad \mbox{for all}\ z\in C(\bH).
\]
This shows that $E_\bH$ is $\bhG$-equivariant 
if and only if the image of the Poisson integral of $\bH$ 
consists of scalars. 
By \cite[Corollary 3.9]{I1}, 
it is equivalent to that $\bH$ is of Kac type. 
\end{proof}

\begin{lem}\label{lem: factor}
Let $\bK$ be a quantum subgroup of Kac type in $\bG$.  
Then $\Th\circ E_\bK=\Th$. 
\end{lem}
\begin{proof}
By Lemma \ref{lem: Kac}, $E_\bK$ is $\bhG$-equivariant. 
For any $x\in\lG$, we have 
\begin{align*}
\Th(E_\bK(x))
=&\,
(\id\oti h)(\be(E_\bK(x)))
=
(\id\oti h)((\id\oti E_\bK)(\be(x)))
\\
=&\,
(\id\oti h)(\be(x))
=
\Th(x). 
\end{align*}
\end{proof}

We prepare the notion of maximality for a quantum subgroup 
of Kac type as follows. 

\begin{defn}
Let $\bG$ be a compact quantum group. 
We say that 
a quantum subgroup $\bH$ of Kac type is \textit{maximal} 
if for any quantum subgroup $\bK$ of Kac type, 
we have $L^\infty(\bH\setm \bG)\subs L^\infty(\bK\setm\bG)$. 
\end{defn}

\begin{lem}\label{lem: unique-max}
Let $\bG$ be a compact quantum subgroup. 
If there exists a maximal quantum subgroup of Kac type, 
it is unique in the following sense. 
Let $\{\bH,r_\bH\}$ and 
$\{\bK,r_\bK\}$ be maximal quantum subgroups of Kac type. 
Then there exists a $*$-isomorphism $\th\col \CH\ra \CK$ 
such that $r_\bK=\th\circ r_\bH$. 
\end{lem}
\begin{proof}
By definition, 
we have $L^\infty(\bH\setm \bG)=L^\infty(\bK\setm\bG)$. 
Then by Lemma \ref{lem: subs}, we have the desired 
$*$-isomorphism. 
\end{proof}

\begin{thm}\label{thm: poisson}
Let $\bG$ be a co-amenable compact quantum group. 
Assume that its fusion algebra is commutative. 
Then the following statements hold. 
\begin{enumerate}

\item 
There exists a unique maximal quantum subgroup of Kac type $\bH$. 

\item 
The Poisson integral 
$\Th\col L^\infty(\bH\setm\bG)\ra H^\infty(\bhG)$ 
is an isomorphism. 

\end{enumerate}
\end{thm}
\begin{proof}
(1). 
By using Theorem \ref{thm: exco}, Lemma \ref{lem: Kac} 
and Lemma \ref{lem: prepoisson}, 
there exists a quantum subgroup $\bH$ of Kac type 
such that $B=L^\infty(\bH\setm\bG)$. 
We show the maximality of $\bH$. 
Let $\bK$ be another quantum subgroup of Kac type. 
By Lemma \ref{lem: factor}, $\Th\circ E_{\bK}=\Th$. 
Since $E_B=\La\circ\Th$, we have $E_B\circ E_{\bK}=E_B$. 
Let $e_{\bK}$ be the Jones projection associated with $E_\bK$. 
Then $e_B e_\bK=e_B$, and hence $e_\bK e_B=e_B$. 
It yields $B\subs L^\infty(\bK\setm \bG)$. 
Hence $\bH$ is maximal. 

(2). 
It follows from Lemma \ref{lem: prepoisson} (3). 
\end{proof}

\begin{rem}
We have realized the existence of 
the maximal quantum subgroup of Kac type by studying the Poisson integral. 
However, 
the notion has been already introduced in \cite[Appendix A]{So}, 
where it is called the canonical Kac quotient. 
The canonical Kac quotient is the function algebra on 
the maximal quantum subgroup of Kac type. 
The existence is proved 
for an arbitrary compact quantum group \cite[Proposition~A.1]{So}. 
\end{rem}

\subsection{$q$-deformed classical compact Lie groups.}

We determine the maximal quantum subgroup of Kac type 
in a $q$-deformed classical compact Lie group ($0<q<1$). 
In order to do so, we freely make use of terminology and results 
in \cite{KS} such as construction of quantum universal 
enveloping algebras, quantized function algebras and so on. 

Let $\mfg$ be a complex classical simple Lie algebra 
with the simple roots $\{\al_i\}_{i=1}^n$. 
We denote by $\bG$ 
the corresponding classical compact Lie group 
\cite[Definition~1.2.4 of Chapter 3]{KS}. 
Let $U_q(\mfg)$ be 
the quantum universal enveloping algebra associated 
with $\mfg$ \cite[Definition 7.1.1 of Chapter 2]{KS}. 
We equip $U_q(\mfg)$ with the Hopf $*$-algebra structure 
as in \cite[Proposition 1.2.3 of Chapter 3]{KS}. 
We denote by $W$ the Weyl group and by $s_i$ the simple reflection 
with respect to $\al_i$. 

For a dominant integral weight $\la$, 
we denote by $L_q(\la)$ the 
irreducible $U_q(\mfg)$-module with the highest weight $\la$. 
Fix an orthonormal basis $\{\xi_\mu^r\}_r$ of each weight space 
$L_q(\la)_\mu$ with the weight $\mu$. 
We define 
$C_{\xi_\mu^r,\xi_\nu^s}^\la\in \AG\subs U_q(\mfg)^*$ 
by
\[
C_{\xi_\mu^r,\xi_\nu^s}^\la (x)=(x \xi_\nu^s, \xi_\mu^r)
\quad\mbox{for}\ 
x\in U_q(\mfg). 
\]
Then we have 
$A(\bG_q)=\spa\{C_{\xi_\mu^r,\xi_\nu^s}^\la\}_{\la,\mu,\nu,r,s}$. 
On the involution, we have 
\begin{equation}\label{eq: involution}
(C_{\xi_\mu^r,\xi_\nu^s}^\la)^*
=q^{(\mu-\nu,\vrho)}C_{\xi_{-\mu}^r,\xi_{-\nu}^s}^{-w_0\la}, 
\end{equation}
where $\vrho=(1/2)\sum_{i=1}^n \al_i$ is the Weyl vector 
and $w_0\in W$ is the longest element. 
Let $\ka$ and $\{\ta_t\}_{t\in\R}$ be the antipode 
and the scaling automorphism group of $\bG_q$ 
introduced in \S \ref{subsect: cqg}. 
From (\ref{eq: antipode}), $\ka^2=\ta_{-i}$ 
and (\ref{eq: involution}), 
we obtain 
\[
\ta_t(C_{\xi_\mu^r,\xi_\nu^s}^\la)
=q^{i(\mu-\nu,2\vrho)t}C_{\xi_\mu^r,\xi_\nu^s}^\la.
\] 
We denote by $C(\bG_q)$ 
the $C^*$-completion of $A(\bG_q)$ with respect to 
the universal norm. 
Set $q_i=q^{(\al_i,\al_i)/2}$, $i=1,\dots,n$. 
Let $U_{q_i}(su(2))$ be the quantized universal enveloping algebra 
of the Lie algebra $su(2)$ with the deformation parameter $q_i$. 
The canonical embedding $U_{q_i}(su(2))\ra U_q(\mfg)$ 
induces the restriction map 
$r_i\col A(\bG_q) \ra A(SU_{q_i}(2))$ 
\cite[Subsection 6.1 of Chapter 3]{KS}. 
Since $C(\bG_q)$ is a universal $C^*$-algebra, 
$r_i$ extends to the $*$-homomorphism $C(\bG_q)\ra C(SU_{q_i}(2))$. 
Take a canonical infinite dimensional 
irreducible representation 
$\pi_i\col C(SU_{q_i}(2))\ra B(\el_2)$ 
defined in \cite[Proposition 4.1.1 of Chapter 3]{KS}. 
We note that 
the counit $\vep_i$ of $C(SU_{q_i}(2))$ factors through 
$\Ima (\pi_i)$, that is, there exists a $*$-homomorphism 
$\eta_i\col \Ima (\pi_i)\ra\C$ such that $\eta_i\circ\pi_i=\vep_i$. 
Indeed, let $p\col B(\el_2)\ra \meQ$ be the canonical 
surjection onto the Calkin algebra $\meQ$. 
Let $S\in B(\el_2)$ be the unilateral shift. 
By definition of $\pi_i$, we see that 
$p(\Ima (\pi_i))$ is a commutative $C^*$-algebra 
generated by a unitary $p(S)$. 
Take the character $\om\col p(\Ima (\pi_i))\ra \C$ defined 
by $\om(p(S))=1$. 
Then the character 
$\eta_i:=\om\circ p\col \Ima(\pi_i)\ra\C$ has the desired property. 

Let $\{\T,r_\T\}$ be the maximal torus subgroup of $\bG_q$, 
where $r_\T\col \CGq\ra C(\T)$ is a restriction map. 
Any one dimensional $*$-homomorphism of $\CGq$ is given by 
$\pi_t:=\chi_t\circ r_\T$, where $\chi_t\col C(\T)\ra\C$ 
is the evaluation at $t\in\T$. 
Then for an element $w\in W\setm\{e\}$ with a reduced decomposition 
$w=s_{i_1}\cdots s_{i_k}$ and $t\in\T$, 
we define 
\[
\pi_{w,t}
=(\pi_{i_1}\circ r_{i_1}\oti\dots\oti \pi_{i_k}\circ r_{i_k}\oti\pi_t)
\circ \de_{\bG_q}^{(k)}, 
\]
where $\de_{\bG_q}^{(k)}\col \CGq\ra \CGq^{\oti (k+1)}$ is recursively defined 
by 
$\de_{\bG_q}^{(k)}=(\de_{\bG_q}\oti\id^{\oti(k-1)})\circ\de_{\bG_q}^{(k-1)}$ 
and $\de_{\bG_q}^{(1)}=\de_{\bG_q}$. 
Then $\pi_{w,t}$ is an irreducible $*$-homomorphism of $\CGq$ 
and does not depend on the choice of the reduced decomposition of $w$ 
up to equivalence \cite[Theorem 6.2.1 of Chapter 3]{KS}. 
Every irreducible $*$-homomorphism of $\CGq$ 
is equivalent to some $\pi_{w,t}$ \cite[Theorem 6.2.7 of Chapter 3]{KS}. 

\begin{lem}
Let $\bG_q$ be the $q$-deformation of a classical compact Lie group $\bG$. 
Then its maximal quantum subgroup of Kac type is the maximal torus $\T$. 
\end{lem}
\begin{proof}
Let $\{\bH,r_\bH\}$ be a quantum subgroup of Kac type. 
We will show that $\bH$ is a subgroup of $\T$. 
Assume that $\CH$ is not commutative. 
Then there exists an irreducible $*$-homomorphism 
$\rho\col \CH\ra B(H_\rho)$ with $\dim H_\rho\geq2$. 
Set $\pi=\rho\circ r_\bH$. 
We may assume $\pi=\pi_{w,t}$
for some element $w\in W\setm\{e\}$ with a reduced expression 
$w=s_{i_1}\cdots s_{i_k}$ and $t\in \T$. 
Then we consider the irreducible $*$-homomorphism 
$(\id\oti \eta_{i_2}\cdots \oti \eta_{i_k})\circ \pi
=(\pi_{i_1}\circ r_{i_1}\oti \pi_t)\circ\de_{\bG_q}=\pi_{s_{i_1},t}$, 
which factors through $\CH$. 
Hence we may assume $\pi=\rho\circ r_\bH=\pi_{s_j,t}$ 
for some $1\leq j\leq n$ and $t\in \T$. 
Since the scaling automorphism of $\CH$ is trivial, 
we have $r_\bH\circ \ta_t=r_\bH$ for all $t\in\R$ 
by Lemma \ref{lem: algsub-scaling}. 
Let $\om_j$ be the fundamental weight for $\al_j$. 
Setting $\la=\mu=\om_j$ and $\nu=\om_j-\al_j$, 
we have 
$\ta_t(C_{\xi_{\om_j},\xi_{\om_j-\al_j}}^{\om_j})
=q_j^{2it}C_{\xi_{\om_j},\xi_{\om_j-\al_j}}^{\om_j}$. 
Hence $r_\bH(C_{\xi_{\om_j},\xi_{\om_j-\al_j}}^{\om_j})=0$. 
In particular, $\pi_{s_j,t}(C_{\xi_{\om_j},\xi_{\om_j-\al_j}}^{\om_j})=0$. 
This is, however, in contradiction 
with $\pi_{s_j,t}(C_{\xi_{\om_j},\xi_{\om_j-\al_j}}^{\om_j})\neq0$ 
which is shown by direct computation. 

Therefore $\CH$ is a commutative $C^*$-algebra, 
that is, $\bH$ is an ordinary compact group. 
For $x\in\bH$, we write $\chi_x^\bH$ 
for the character $*$-homomorphism. 
For $t\in\T$, we also use $\chi_t^\T$ as well. 
By the above arguments, for any $x\in \bH$ 
there exists an element $t(x)\in\T$ such that 
$\chi_x^\bH\circ r_\bH=\chi_{t(x)}^\T\circ r_\T$. 
This shows that the map $t\col \bH\ra \T$ is continuous. 
For any $x, y\in\bH$, we have 
\begin{align*}
\chi_{t(xy)}^\T\circ r_\T
=&\,
\chi_{xy}^\bH\circ r_\bH
=
(\chi_x^\bH\oti \chi_y^\bH)\circ\de_\bH\circ r_\bH
=
(\chi_x^\bH\oti \chi_y^\bH)\circ(r_\bH\oti r_\bH)
\circ \de_{\bG_q}
\\
=&\,
(\chi_{t(x)}^\T \oti\chi_{t(y)}^\T)\circ (r_\T\oti r_\T)
\circ \de_{\bG_q}
=
(\chi_{t(x)}^\T \oti\chi_{t(y)}^\T)\circ \de_\T \circ r_\T
\\
=&\,
\chi_{t(x)t(y)}^\T\circ r_\T. 
\end{align*}
Hence the map $t\col \bH\ra \T$ is a group homomorphism. 
Next we show that $t$ is injective. 
For $x\in \bH$, 
suppose that $\chi_x^\bH\circ r_\bH=\chi_e^\T\circ r_\T$ holds. 
Since it is equal to the counit $\vep$ of $\CGq$ 
by Lemma \ref{lem: algsub-scaling}, we have 
\[
\chi_x^\bH\circ r_\bH
=
\chi_e^\T\circ r_\T
=
\vep 
=
\chi_e^\bH\circ r_\bH. 
\]
Hence $x=e$, and the map $t$ is injective. 
Since the left action of $x\in \bH$ is given by 
$
(\chi_x^\bH \circ r_\bH\oti\id)\circ\de_{\bG_q}
=(\chi_{t(x)}^\T\circ r_\T\oti\id)\circ \de_{\bG_q}
$, 
we have $L^\infty(\T\setm\bG_q)\subs L^\infty(\bH\setm \bG_q)$. 
Therefore $\T$ is maximal. 
\end{proof}

By Theorem \ref{thm: poisson} and the previous lemma, 
we obtain the following corollary. 
\begin{cor}\label{cor: flag}
Let $\bG_q$ be the $q$-deformation 
of a classical compact Lie group $\bG$. 
Then the Poisson integral 
$\Th\col L^\infty(\T\setm\bG_q)\ra H^\infty(\widehat{\bG_q})$ 
is an isomorphism. 
\end{cor}

\textbf{Acknowledgments.} 
The author is grateful to Yasuyuki Kawahigashi and Masaki Izumi 
for encouragement and various useful comments. 
He would also like to thank Stefaan Vaes for informing him 
of canonical Kac quotients introduced in \cite{So}.

\end{document}